\documentclass{amsart}

 \newtheorem{thm}{Theorem}[section]

 \newtheorem{prop}[thm]{Proposition}
 \theoremstyle{definition}
 
 \theoremstyle{remark}

 \numberwithin{equation}{section}

\newcommand{\bpr}{\begin{proof}[Proof]}  
\newcommand{\epr}{\end{proof}}
\newcommand{\beq}{\begin{equation}}
\newcommand{\eeq}{\end{equation}}
\newcommand{\bce}{\begin{center}}
\newcommand{\ece}{\end{center}}
\newcommand{\be}{\begin{enumerate}}  
\newcommand{\ee}{\end{enumerate}}
\newcommand{\difft}{\frac{d}{dt}}

\newcommand{\vectfuenf}[5]{\left[\begin{matrix}#1\\  #2\\ #3\\ #4\\ #5\end{matrix}\right]}

\DeclareMathOperator*{\dist}{dist}
\DeclareMathOperator*{\diver}{div}

\def\al{\alpha}
\def\pa{\partial}
\def\om{\omega}
\def\si{\sigma}

\def\la{\lambda}

\def\ep{\varepsilon}
\def\de{\delta}
\def\ph{\varphi}

\def\ka{\kappa}
\def\ga{\gamma}
\def\th{\theta}
\def\De{\Delta}

\def\Om{\Omega}
\def\vth{\vartheta}
\def\R{\mathbb R}

\def\N{\mathbb N}

\def\X{\mathbb X}
\def\B{\mathbb B}
\def\E{\mathbb E}
\def\L{\mathbb L}

\def\calA{\mathcal A}

\def\calT{\mathcal T}

\def\calM{\mathcal M}

\def\nn{\nonumber}

\begin{document}
%
%
%
%
%
%
%
%
%
\title[Conserved Penrose-Fife type models]
 {On conserved Penrose-Fife type models}
\author[J. Pr\"{u}ss]{Jan Pr\"{u}ss}

\address{
Theodor-Lieser-Str. 5\\
06120 Halle\\
Deutschland}

\email{jan.pruess@mathematik.uni-halle.de}

\author[M. Wilke]{Mathias Wilke}
\address{Theodor-Lieser-Str. 5\\
06120 Halle\\
Deutschland}

\email{mathias.wilke@mathematik.uni-halle.de (Corresponding author)}


\date{\today}
\dedicatory{Dedicated to Herbert Amann on the occasion of his 70th birthday}

\begin{abstract}
In this paper we investigate quasilinear parabolic systems of conserved Penrose-Fife type. We show maximal $L_p$ - regularity for this problem
with inhomogeneous boundary data. Furthermore we prove global existence of a solution, provided that the absolute temperature is bounded from below and above. Moreover, we apply the Lojasiewicz-Simon inequality to establish the convergence of solutions to a steady state as time tends to infinity.
\end{abstract}

\maketitle

\section{Introduction and the Model}
We are interested in the conserved Penrose-Fife type
equations
\beq\begin{split}\label{PF0}\pa_t\psi=\De\mu,\quad \mu=-\De\psi+\Phi'(\psi)-\la'(\psi)\vth,&\quad t\in J,\ x\in\Om,\\
\pa_t\left(b(\vth)+\la(\psi)\right)-\De\vth=0,&\quad t\in J,\
x\in\Om,
\end{split}\eeq
where $\vartheta=1/\theta$ and $\theta$ denotes the absolute temperature of the system, $\psi$ is the order parameter and $\Omega\subset\mathbb{R}^n$ is a bounded domain with boundary $\partial\Omega\in C^4$. The function $\Phi'$ is the derivative of the physical potential, which characterizes the different phases of the system. A typical example is the \emph{double well} potential $\Phi(s)=(s^2-1)^2$ with the two distinct minima $s=\pm1$. Typically, the nonlinear function $\lambda$ is a polynomial of second order.

For an explanation of \eqref{PF0} we will follow the lines of \textsc{Alt \& Pawlow} \cite{AltPaw} (see also \textsc{Brokate} \& \textsc{Sprekels} \cite[Section 4.4]{BroSpr}). We start with the rescaled Landau-Ginzburg functional (total Helmholtz free energy)
$$\mathcal{F}(\psi,\theta)=\int_\Om\left(\frac{\gamma(\theta)}{2\theta}|\nabla\psi|^2+\frac{f(\psi,\theta)}{\theta}\right)\
dx,$$
where the free energy density $F(\psi,\theta):=\frac{\gamma(\theta)}{2}|\nabla\psi|^2+f(\psi,\theta)$ is rescaled by $1/\theta$. The reduced chemical potential $\mu$ is given by the variational derivative of $\mathcal{F}$ with respect to $\psi$, i.e.
$$\mu=\frac{\de \mathcal{F}}{\de\psi}(\psi,\theta)=\frac{1}{\theta}\left(-\gamma(\theta)\De\psi+\frac{\partial f(\psi,\theta)}{\partial \psi}\right).$$
Assuming that $\psi$ is a conserved quantity, we have the conservation law
    $$\partial_t\psi+{\rm div} j=0.$$
Here $j$ is the flux of the order parameter $\psi$, for which we choose the well accepted constitutive law $j=-\nabla\mu$, i.e.\ the phase transition is driven by the chemical potential $\mu$ (see \cite[(4.4)]{BroSpr}). The kinetic equation for $\psi$ thus reads
    $$\partial_t\psi=\Delta\mu,\quad \mu=\frac{1}{\theta}\left(-\gamma(\theta)\De\psi+\frac{\partial f(\psi,\theta)}{\partial \psi}\right).$$
If the volume of the system is preserved, the internal energy $e$ is given by the variational derivative
    $$e=\frac{\delta\mathcal{F}(\psi,\theta)}{\delta (1/\theta)}.$$
This yields the expression
    $$e(\psi,\theta)=
    f(\psi,\theta)-\theta\frac{\partial f(\psi,\theta)}{\partial\theta}+\frac{1}{2}\left(\gamma(\theta)-
    \theta\frac{\partial\gamma(\theta)}{\partial\theta}\right)|\nabla\psi|^2.$$
It can be readily checked that the \textsc{Gibbs} relation
    $$e(\psi,\theta)=F(\psi,\theta)-\theta\frac{\partial F(\psi,\theta)}{\partial\theta}.$$
holds. If we assume that no mechanical stresses are active, the internal energy $e$ satisfies the conservation law
    $$\partial_t e+{\rm div}q=0,$$
where $q$ denotes the heat flux of the system. Following \textsc{Alt \& Pawlow} \cite{AltPaw}, we assume that $q=\nabla\left(\frac{1}{\theta}\right)$, so that the kinetic equation for $e$ reads
    $$\partial_t e+\Delta\left(\frac{1}{\theta}\right)=0.$$
Let us now assume that $\gamma(\theta)=\theta$ and $f(\psi,\theta)=\theta\Phi(\psi)-\lambda(\psi)-\theta\log\theta$. In this case we obtain
$e=\theta-\lambda(\psi)$ and
    $$\mu=-\Delta\psi+\Phi'(\psi)-\lambda'(\psi)\frac{1}{\theta},$$
hence system \eqref{PF0} for $\vartheta=1/\theta$ and $b(s)=-1/s$, $s>0$. Suppose $(j|\nu)=(q|\nu)=0$ on $\partial\Omega$ with $\nu=\nu(x)$ being the outer unit normal in $x\in\partial\Omega$. This yields the boundary conditions
$\pa_\nu\mu=0$ and $\pa_\nu\vth=0$ for the chemical potential $\mu$
and the function $\vth$, respectively. Since \eqref{PF0}
is of fourth order with respect to the function $\psi$ we need an
additional boundary condition. An appropriate and classical one from
a variational point of view is $\pa_\nu\psi=0$. Finally, this yields
the initial-boundary value problem
\beq\label{PF}\begin{split}\pa_t\psi-\De\mu=f_1,\quad \mu=-\De\psi+\Phi'(\psi)-\la'(\psi)\vth,&\quad t\in J,\ x\in\Om,\\
\pa_t\left(b(\vth)+\la(\psi)\right)-\De\vth=f_2,&\quad t\in J,\ x\in\Om,\\
\pa_\nu\mu=g_1,\ \pa_\nu\psi=g_2,\ \pa_\nu\vth=g_3,&\quad t\in J,\ x\in\pa\Om,\\
\psi(0)=\psi_0,\ \vth(0)=\vth_0,&\quad t=0,\ x\in\Om,
\end{split}\eeq
The functions $f_j, g_j, \psi_0, \vth_0, \Phi, \la$ and $b$ are
given. Note that if $\theta$ has only a small deviation from a constant value $\theta_*>0$, then the term $1/\theta$ can be linearized around $\theta_*$ and \eqref{PF} turns into the nonisothermal Cahn-Hilliard equation for the order parameter $\psi$ and the relative temperature $\theta-\theta_*$, provided $b(s)=-1/s$.

In the case of the Penrose-Fife equations, \textsc{Brokate \& Sprekels}
\cite{BroSpr} and \textsc{Zheng} \cite{Zhe2} proved global
well-posedness in an $L_2$-setting if the spatial dimension is
equal to 1. \textsc{Sprekels \& Zheng} showed global well-posedness
of the non-conserved equations (that is $\pa_t\psi=-\mu$) in
higher space dimensions in \cite{SprZhe}, a similar result can be
found in the article of \textsc{Laurencot} \cite{Lau}. Concerning
asymptotic behavior we refer to the articles of \textsc{Kubo, Ito
\& Kenmochi} \cite{KubItoKen}, \textsc{Shen \& Zheng}
\cite{SheZhe}, \textsc{Feireisl \& Schimperna} \cite{FeiSchim} and
\textsc{Rocca \& Schimperna} \cite{Rocca1}. The last two authors
studied well-posedness and qualitative behavior of solutions to
the non-conserved Penrose-Fife equations. To be precise, they
proved that each solution converges to a steady state, as time
tends to infinity. \textsc{Shen \& Zheng} \cite{SheZhe}
established the existence of attractors for the non-conserved
equations, whereas \textsc{Kubo, Ito \& Kenmochi} \cite{KubItoKen}
studied the non-conserved as well as the conserved Penrose-Fife
equations. Beside the proof of global well-posedness in the sense
of weak solutions they also showed the existence of a global
attractor. Finally, we want to mention that the physical potential $\Phi$ may also be of logarithmic type, such that $\Phi'(s)$ has singularities at $s=\pm1$. This forces the order parameter to stay in the physically reasonable interval $(-1,1)$, provided that the initial value $\psi(0)=\psi_0\in (-1,1)$. In general, such a result cannot be obtained in the case of the double well potential, since there is no maximum principle available for the fourth order equation $\eqref{PF}_1$. For a result on global existence, uniqueness and asymptotic behaviour of solutions of the \emph{Cahn-Hilliard equation} in case of a logarithmic potential, we refer the reader to \textsc{Abels} \& \textsc{Wilke} \cite{AbWi}. However, in this paper we will only deal with smooth potentials.

In the following sections we will prove well-posedness of
\eqref{PF} for solutions in the maximal $L_p$-regularity classes
$$\psi\in H_p^1(J;L_p(\Om))\cap L_p(J;H_p^4(\Om)),$$
$$\vth\in H_p^1(J;L_p(\Om))\cap L_p(J;H_p^2(\Om)),$$
where $J=[0,T]$, $T>0$.
In Section 2 we investigate a linearized version of \eqref{PF} and prove maximal $L_p$-regularity. Section 3 is devoted to local well-posedness of \eqref{PF}. To this end we apply the contraction mapping principle. In Section 4, we show that the solution exists globally in time, provided that the absolute temperature $\vth$ is uniformly bounded from below and above. Finally, in Section 5, we study the asymptotic behavior of the solution to \eqref{PF} as $t\to\infty$. The Lojasiewicz-Simon inequality will play an important role in the analysis.

\section{The Linear Problem}\label{LinProblem}

In this section we deal with a linearized version of \eqref{PF}.
\beq\begin{split}\label{linPF}
\pa_tu+\De^2u+\De(\eta_1 v)=f_1,&\quad t\in J,\ x\in\Om,\\
\pa_tv-a_0\De v+\eta_2\pa_tu=f_2,&\quad t\in J,\ x\in\Om,\\
\pa_\nu\De u+\pa_\nu(\eta_1 v)=g_1,&\quad t\in J,\ x\in\pa\Om,\\
\pa_\nu u=g_2,\ \pa_\nu v=g_3,&\quad t\in J,\ x\in\pa\Om,\\
u(0)=u_0,\ v(0)=v_0,&\quad t=0,\ x\in\Om.\end{split}\eeq Here
$\eta_1=\eta_1(x),\eta_2=\eta_2(x),a_0=a_0(x)$ are given functions
such that \beq\label{condeta}\eta_1\in B_{pp}^{4-4/p}(\Omega),\ \eta_2\in B_{pp}^{2-2/p}(\Omega)\quad\text{and}\quad a_0\in C(\overline{\Om}).\eeq We
assume furthermore that $a_0(x)\ge \si>0$ for all
$x\in\overline{\Om}$ and some constant $\si>0$. Hence equation
$\eqref{linPF}_2$ does not degenerate. We are interested in
solutions
$$u\in H_p^1(J;L_p(\Om))\cap L_p(J;H_p^4(\Om))=:E_1(T)$$
and
$$v\in H_p^1(J;L_p(\Om))\cap L_p(J;H_p^2(\Om))=:E_2(T)$$
of \eqref{linPF}. By the well-known trace theorems (cf.
\cite[Theorem 4.10.2]{Ama}) \beq\label{trace}E_1(T)\hookrightarrow
C(J;B_{pp}^{4-4/p}(\Om))\quad\text{and}\quad E_2(T)\hookrightarrow
C(J;B_{pp}^{2-2/p}(\Om)),\eeq we necessarily have $u_0\in
B_{pp}^{4-4/p}(\Om)=:X_{\gamma}^1$, $v_0\in B_{pp}^{2-2/p}(\Om)=:X_\gamma^2$
and the compatibility conditions
$$\pa_\nu\De u_0+\pa_\nu(\eta_1 v_0)=g_1|_{t=0},\quad \pa_\nu u_0=g_2|_{t=0},\quad\text{as well
as}\quad \pa_\nu v_0=g_3|_{t=0},$$ whenever $p>5$, $p>5/3$ and $p>3$, respectively (cf. \cite[Theorem 2.1]{DHP07}). In the sequel we will assume that $p>(n+2)/2$ and $p\ge 2$. This yields the embeddings
    $$B_{pp}^{4-4/p}(\Omega)\hookrightarrow H_p^2(\Omega)\cap C^1(\bar{\Omega})\ \mbox{and}\ B_{pp}^{2-2/p}(\Omega)\hookrightarrow H_p^1(\Omega)\cap C(\bar{\Omega}).$$
We are going to prove the following theorem.
\begin{thm}\label{linthm}
Let $n\in\N$, $\Omega\subset\R^n$ a bounded domain with boundary $\partial\Omega\in C^4$ and let $p>(n+2)/2$, $p\ge 2$, $p\neq 3,5$. Assume in addition that $\eta_1\in B_{pp}^{4-4/p}(\Omega)$, $\eta_2\in B_{pp}^{2-2/p}(\Omega)$ and $a_0\in C(\bar{\Omega})$, $a_0(x)\ge \sigma>0$ for all $x\in\bar{\Omega}$. Then the linear problem \eqref{linPF} admits a unique solution
    $$(u,v)\in H_p^1(J_0;L_p(\Omega)^2)\cap L_p(J_0;(H_p^4(\Omega)\times H_p^2(\Omega))),$$
if and only if the data are subject to the following conditions.
    \begin{enumerate}
    \item $f_1,f_2\in L_p(J_0;L_p(\Omega))=X(J_0)$,
    \item $g_1\in W_p^{1/4-1/4p}(J_0;L_p(\partial\Omega))\cap L_p(J_0;W_p^{1-1/p}(\partial\Omega))=Y_1(J_0)$,
    \item $g_2\in W_p^{3/4-1/4p}(J_0;L_p(\partial\Omega))\cap L_p(J_0;W_p^{3-1/p}(\partial\Omega))=Y_2(J_0)$,
    \item $g_3\in W_p^{1/2-1/2p}(J_0;L_p(\partial\Omega))\cap L_p(J_0;W_p^{1-1/p}(\partial\Omega))=Y_3(J_0)$,
    \item $u_0\in B_{pp}^{4-4/p}(\Omega)=X_{\gamma}^1$, $v_0\in B_{pp}^{2-2/p}(\Omega)=X_\gamma^2$,
    \item $\pa_\nu\De u_0+\pa_\nu(\eta_1 v_0)=g_1|_{t=0},\ p>5$,
    \item $\pa_\nu u_0=g_2|_{t=0}$, $\pa_\nu v_0=g_3|_{t=0},\ p>3$.
    \end{enumerate}
\end{thm}
\begin{proof}
Suppose that the function $u\in E_1(T)$ in \eqref{linPF} is already
known. Then in a first step we will solve the linear heat equation
\beq\label{heateq} \pa_tv-a_0\De v=f_2-\eta_2\pa_tu,\eeq subject to
the boundary and initial conditions $\pa_\nu v=g_3$ and $v(0)=v_0$.
By the properties of the function $a_0$ we may apply \cite[Theorem 2.1]{DHP07} to obtain a unique solution $v\in E_2(T)$ of \eqref{heateq},
provided that $f_2\in L_p(J\times\Om)$, $v_0\in
B_{pp}^{2-2/p}(\Om)$,
$$g_3\in W_p^{1/2-1/2p}(J;L_p(\partial\Omega))\cap L_p(J;W_p^{1-1/p}(\partial\Omega))=:Y_3(J),$$ and
the compatibility condition $\pa_\nu v_0=g_3|_{t=0}$ if $p>3$ is
valid. The solution may then be represented by the variation of
parameters formula
\begin{align}\label{solvth}
v(t)&=v_1(t)-\int_0^te^{-A(t-s)}\eta_2\pa_tu(s)\ ds,
\end{align}
where $A$ denotes the $L_p$-realization of the differential
operator $\calA(x)=-a_0(x)\De_N$, $\De_N$ means the
Neumann-Laplacian and $e^{-At}$ stands for the bounded analytic
semigroup, which is generated by $-A$ in $L_p(\Om)$. Furthermore
the function $v_1\in E_2(T)$ solves the linear problem
$$\pa_tv_1-a_0\De v_1=f_2,\quad \pa_\nu v_1=g_3,\quad v_1(0)=v_0.$$
We fix a function $w^*\in E_1(T)$ such that $w^*|_{t=0}=u_0$ and make
use of \eqref{solvth} and the fact that $(u-w^*)|_{t=0}=0$ to
obtain
\begin{align*}
v(t)=v_1(t)+v_2(t)-(\pa_t+A)^{-1}\eta_2\pa_t(u-w^*)
\end{align*}
with $v_2(t):=-\int_0^te^{-A(t-s)}\eta_2\pa_t w^*$. Set
$v^*=v_1+v_2\in E_2(T)$ and
$$F(u)=-(\pa_t+A)^{-1}\eta_2\pa_t(u-w^*).$$
Then we may reduce \eqref{linPF} to the problem
\beq\begin{split}\label{linu}
\pa_tu+\De^2u=\De G(u)+f_1,&\quad t\in J,\ x\in\Om,\\
\pa_\nu\De u=\pa_\nu G(u)+g_1,&\quad t\in J,\ x\in\pa\Om,\\
\pa_\nu u=g_2&\quad t\in J,\
x\in\pa\Om,\\
u(0)=u_0,&\quad t=0,\ x\in\Om,
\end{split}\eeq
where $G(u):=-\eta_1(F(u)+v^*)$. For a given $T\in (0,T_0]$ we set
$$_0E_1(T)=\{u\in E_1(T):u|_{t=0}=0\}$$
and
$$E_0(T):=X(T)\times Y_1(T)\times Y_2(T)$$
$$_0E_0(T):=\{(f,g,h)\in E_0(T):g|_{t=0}=h|_{t=0}=0\},$$
where $X(T):=L_p((0,T)\times\Om)$,
$$Y_1(T):=W_p^{1/4-1/4p}(0,T;L_p(\partial\Omega))\cap L_p(0,T;W_p^{1-1/p}(\partial\Omega)),$$
and
$$Y_2(T):=W_p^{3/4-1/4p}(0,T;L_p(\partial\Omega))\cap L_p(0,T;W_p^{3-1/p}(\partial\Omega)).$$
The spaces $E_1(T)$ and $E_0(T)$ are endowed with the canonical norms
$|\cdot|_1$ and $|\cdot|_0$, respectively. We introduce the new function $\tilde{u}:=u-w^*\in\! _0E_1(T)$ and we set
    $$\tilde{F}(\tilde{u}):=-(\pa_t+A)^{-1}\eta_2\pa_t\tilde{u}$$
as well as $\tilde{G}(\tilde{u}):=-\eta_1\tilde{F}(\tilde{u})$. If $u\in E_1(T)$ is a solution of \eqref{linu}, then the function $\tilde{u}\in\! _0E_1(T)$ solves the problem
\beq\begin{split}\label{linutilde}
\pa_t\tilde{u}+\De^2\tilde{u}=\De \tilde{G}(\tilde{u})+\tilde{f}_1,&\quad t\in J,\ x\in\Om,\\
\pa_\nu\De \tilde{u}=\pa_\nu \tilde{G}(\tilde{u})+\tilde{g}_1,&\quad t\in J,\ x\in\pa\Om,\\
\pa_\nu \tilde{u}=\tilde{g}_2&\quad t\in J,\
x\in\pa\Om,\\
\tilde{u}(0)=0,&\quad t=0,\ x\in\Om,
\end{split}\eeq
with the modified data
    $$\tilde{f}_1:=f_1-\Delta(\eta_1 v^*)-\partial_t w^*-\Delta^2 w^*\in X(T),$$
    $$\tilde{g}_1:=g_1-\partial_\nu (\eta v^*)-\partial_\nu\Delta w^*\in\! _0Y_1(T),$$
and
    $$\tilde{g}_2:=g_2-\partial_\nu w^*\in\! _0Y_2(T).$$
Observe that by construction we have $\tilde{g}_1|_{t=0}=0$ and $\tilde{g}_2|_{t=0}=0$ if $p>5$ and $p>5/3$, respectively.

Let us estimate the term $\Delta \tilde{G}(u)$ in $L_p(J;L_p(\Omega))$, where $u\in\! _0E_1(T)$. We compute
    \begin{multline*}
    |\Delta\tilde{G}(u)|_{L_p(J;L_p(\Omega))}\le |\tilde{F}(u)\Delta\eta_1|_{L_p(J;L_p(\Omega))}\\
    +2|(\nabla \tilde{F}(u)|\nabla\eta_1)|_{L_p(J;L_p(\Omega))}+|\eta_1\Delta \tilde{F}(u)|_{L_p(J;L_p(\Omega))}.
    \end{multline*}
Since $\eta_1\in B_{pp}^{4-4/p}(\Omega)$ does not depend on the variable $t$, we obtain
    $$|\tilde{F}(u)\Delta\eta_1|_{L_p(J;L_p(\Omega))}\le |\Delta\eta_1|_{L_p(\Omega)}|\tilde{F}(u)|_{L_p(J;L_\infty(\Omega))},$$
    $$|(\nabla \tilde{F}(u)|\nabla\eta_1)|_{L_p(J;L_p(\Omega))}\le |\nabla\eta_1|_{L_\infty(\Omega)}|\nabla\tilde{F}(u)|_{L_p(J;L_p(\Omega))},$$
and
    $$|\eta_1\Delta \tilde{F}(u)|_{L_p(J;L_p(\Omega))}\le |\eta_1|_{L_\infty(\Omega)}|\Delta \tilde{F}(u)|_{L_p(J;L_p(\Omega))}.$$
Therefore we have to estimate $\tilde{F}(u)$ for each $u\in\! _0E_1(T)$ in the topology of the spaces $L_p(J;L_\infty(\Omega))$ and $L_p(J;H_p^2(\Omega))$. Let $u\in\! _0E_1$ and recall that $\tilde{F}(u)$ is defined by $\tilde{F}(u)=-(\pa_t+A)^{-1}\eta_2 \pa_t u$. The operator $(\partial_t+A)^{-1}$ is a bounded linear operator from $L_p(J;L_p(\Omega))$ to $_0H_p^1(J;L_p(\Omega))\cap L_p(J;H_p^2(\Omega))=\! _0E_2(T)$. Moreover, by the trace theorem and by Sobolev embedding, it holds that
    $$_0H_p^1(J;L_p(\Omega))\cap L_p(J;H_p^2(\Omega))\hookrightarrow C(J;B_{pp}^{2-2/p}(\Omega))\hookrightarrow C(J;C(\bar{\Omega})).$$
Note that the bound of $(\partial_t+A)^{-1}$ as well as the embedding constant do not depend on the length of the interval $J=[0,T]\subset [0,T_0]=J_0$, since the time trace at $t=0$ vanishes. With these facts, we obtain
    \begin{align*}
    |(\pa_t+A)^{-1}\eta_2 \partial_t u|_{L_p(J;L_\infty(\Omega))}&\le T^{1/p}|(\pa_t+A)^{-1}\eta_2 \pa_tu|_{L_\infty(J;L_\infty(\Omega))}\\
    &\le T^{1/p}C|(\pa_t+A)^{-1}\eta_2 \pa_tu|_{E_2(T)}\\
    &\le T^{1/p}C|\eta_2\partial_tu|_{L_p(J;L_p(\Omega))}\\
    &\le T^{1/p}C|\eta_2|_{L_{\infty}(\Omega)}|u|_{E_1(T)}.
    \end{align*}
To estimate $\tilde{F}(u)$ in $L_p(J;H_p^2(\Omega))$ we need another representation of $\tilde{F}(u)$. To be precise, we rewrite $\tilde{F}(u)$ as follows
    $$\tilde{F}(u)=-(\pa_t+A)^{-1}\eta_2 \pa_tu=-\partial_t^{1/2}(\pa_t+A)^{-1}\pa_t^{1/2}(\eta_2 u).$$
This is possible, since $u\in\! _0E_1(T)$. Now observe that for each $u\in\! _0E_1$ it holds that $\eta_2 u\in\! _0H_p^{3/4}(J;H_p^1(\Omega))$. This can be seen as follows. First of all, it suffices to show that $\eta_2 u\in L_p(J;H_p^1(\Omega))$, since $\eta_2$ does not depend on the variable $t$. But
    \begin{align*}
    |\eta_2 u|_{L_p(J;H_p^1(\Omega))}&\le |\eta_2\nabla u|_{L_p(J;L_p(\Omega))}+|u\nabla\eta_2|_{L_p(J;L_p(\Omega))}\\
    &\le C\left(|\eta_2|_{L_\infty(\Omega)}|u|_{E_1(T)}+|u|_{L_p(J;L_\infty(\Omega))}|\eta_2|_{H_p^1(\Omega)}\right)\\
    &\le C|u|_{E_1(T)}|\eta_2|_{B_{pp}^{2-2/p}(\Omega)},
    \end{align*}
and this yields the claim, since
    $$u\in\! _0H_p^1(J;L_p(\Omega))\cap L_p(J;H_p^4(\Omega))\hookrightarrow\! _0H_p^{3/4}(J;H_p^1(\Omega)),$$
by the mixed derivative theorem. It follows readily that $\pa_t^{1/2}(\eta_2 u)\in\! _0H_p^{1/4}(J;H_p^1(\Omega))$ and
    $$(\partial_t+A)^{-1}(I+A)^{1/2}\partial_t^{1/2}(\eta_2 u)\in\! _0H_p^{5/4}(J;L_p(\Omega))\cap\! _0H_p^{1/4}(J;H_p^2(\Omega)).$$
Since the operator $(I+A)^{1/2}$ with domain $D((I+A)^{1/2})=H_p^1(\Omega)$ commutes with the operator $(\partial_t+A)^{-1}$,
this yields
    $$(\partial_t+A)^{-1}\partial_t^{1/2}(\eta_2 u)\in\! _0H_p^{5/4}(J;H_p^1(\Omega))\cap\! _0H_p^{1/4}(J;H_p^3(\Omega))$$
for each fixed $u\in\! _0E_1(T)$. By the mixed derivative theorem we obtain furthermore
    $$_0H_p^{5/4}(J;H_p^1(\Omega))\cap\! _0H_p^{1/4}(J;H_p^3(\Omega))\hookrightarrow\! _0H_p^{3/4}(J;H_p^2(\Omega)).$$
Therefore
    $$\tilde{F}(u)=-\partial_t^{1/2}(\partial_t+A)^{-1}\partial_t^{1/2}(\eta_2 u)\in\! _0H_p^{1/4}(J;H_p^2(\Omega)),$$
and there exists a constant $C>0$ being independent of $T>0$ and $u\in\! _0E_1(T)$ such that
    $$|\tilde{F}(u)|_{H_p^{1/4}(J;H_p^2(\Omega))}\le C|u|_{E_1(T)},$$
for each $u\in\! _0E_1(T)$. In particular this yields the estimate
    \begin{align*}
    |\tilde{F}(u)|_{L_p(J;H_p^2(\Omega))}&\le T^{1/2p}|\tilde{F}(u)|_{L_{2p}(J;H_p^2(\Omega))}\\
    &\le T^{1/2p}|\tilde{F}(u)|_{H_p^{1/4}(J;H_p^2(\Omega))}\le T^{1/2p}C|u|_{E_1(T)},
    \end{align*}
by H\"{o}lders inequality and $C>0$ does not depend on the length $T$ of the interval $J$. We have thus shown that
    $$|\Delta\tilde{G}(u)|_{L_p(J;L_p(\Omega))}\le \mu_1(T)C|u|_{E_1(T)},$$
where we have set $\mu_1(T):=T^{1/2p}(1+T^{1/2p})$. Observe that $\mu_1(T)\to 0_+$ as $T\to 0_+$. The next step consists of estimating the term $\partial_\nu \tilde{G}(u)$ in $_0W_p^{1/4-1/4p}(J;L_p(\partial\Omega))\cap L_p(J;W_p^{1-1/p}(\partial\Omega))$. To this end, we recall the trace map
    $$_0H_p^{1/2}(J;L_p(\Omega))\cap L_p(J;H_p^2(\Omega))\hookrightarrow\! _0W_p^{1/4-1/4p}(J;L_p(\partial\Omega))\cap L_p(J;W_p^{1-1/p}(\partial\Omega))$$
for the Neumann derivative on $\partial\Omega$. Therefore, by the results above, it remains to estimate $\tilde{G}(u)$ in $_0H_p^{1/2}(J;L_p(\Omega))$. By the complex interpolation method we have
    $$|w|_{H_p^{1/2}(J;L_p(\Omega))}\le C|w|_{L_p(J;L_p(\Omega))}^{1/2}|w|_{H_p^1(J;L_p(\Omega))}^{1/2}$$
for each $w\in\! _0H_p^1(J;L_p(\Omega))$, and $C>0$ does not depend on $T>0$. Using H\"{o}lders inequality, this yields
    \begin{align*}
    |w|_{H_p^{1/2}(J;L_p(\Omega))}&\le T^{1/4p}C|w|_{L_{2p}(J;L_p(\Omega))}^{1/2}|w|_{H_p^1(J;L_p(\Omega))}^{1/2}\\
    &\le T^{1/4p}C|w|_{H_p^1(J;L_p(\Omega))}.
    \end{align*}
Finally we obtain the estimate
    $$|\tilde{G}(u)|_{H_p^{1/2}(J;L_p(\Omega))}\le T^{1/2p}|\eta_1|_{L_\infty(\Omega)}C|u|_{\E_1(T)},$$
which in turn implies
    \begin{align*}
    |\partial_\nu\tilde{G}(u)|_{Y_1(J)}\le |\tilde{G}(u)|_{H_p^{1/2}(J;L_p(\Omega))}+|\tilde{G(u)}|_{L_p(J;H_p^2(\Omega))}\le \mu_2(T)C|u|_{E_1(T)},
    \end{align*}
where $\mu_2(T):=T^{1/4p}(1+T^{1/4p})$ and $\mu_2(T)\to 0_+$ as $T\to 0_+$. Define two operators $L,B:\!_0E_1(T)\to\! _0E_0(T)$ by means of
    $$Lu:=\begin{bmatrix}\partial_t u+\Delta^2 u\\\partial_\nu\Delta u\\\partial_\nu u\end{bmatrix}\ \mbox{and}\
        Bu:=\begin{bmatrix}\Delta\tilde{G}(u)\\\partial_\nu\tilde{G}(u)\\0\end{bmatrix}.$$
With these definitions, we may rewrite \eqref{linutilde} in the abstract form
    $$Lu=Bu+f,\quad f:=(\tilde{f}_1,\tilde{g}_1,\tilde{g}_2)\in\! _0E_0(T).$$
By \cite[Theorem 2.1]{DHP07}, the operator $L$ is bijective with bounded inverse $L^{-1}$, hence $u\in\! _0E_1(T)$ is a solution of \eqref{linutilde} if and only if $(I-L^{-1}B)u=L^{-1}f$. Observe that $L^{-1}B$ is a bounded linear operator from $_0E_1(T)$ to $_0E_1(T)$ and
    $$|L^{-1}Bu|_{E_1(T)}\le |L^{-1}|_{\mathcal{B}(E_0(T),E_1(T))}|Bu|_{E_0(T)}\le (\mu_1(T)+\mu_2(T))C|u|_{E_1(T)}.$$
Here the constant $C>0$ as well as the bound of $L^{-1}$ are independent of $T>0$. This shows that choosing $T>0$ sufficiently small, we may apply a Neumann series argument to conclude that \eqref{linutilde} has a unique solution $u\in\! _0E_1(T)$ on a possibly small time interval $J=[0,T]$. Since the linear system \eqref{linutilde} is invariant with respect to time shifts, we may set $J=J_0$.
\end{proof}

\section{Local Well-Posedness}\label{LWP}

In this section we will use the following setting. For $T_0>0$, to
be fixed later, and a given $T\in (0,T_0]$ we define
$$\E_1(T):=E_1(T)\times E_2(T),\qquad\hspace{0.05cm}_0\E_1(T):=\{(u,v)\in\E_1(T):(u,v)|_{t=0}=0\}$$
and
$$\E_0(T):=X(T)\times X(T)\times Y_1(T)\times Y_2(T)\times Y_3(T),$$
as well as $$_0\E_0(T):=\{(f_1,f_2,g_1,g_2,g_3)\in\E_0(T):
g_1|_{t=0}=g_2|_{t=0}=g_3|_{t=0}=0\},$$ with canonical norms
$|\cdot|_1$ and $|\cdot|_0$, respectively. The aim of this section
is to find a local solution $(\psi,\vth)\in \E_1(T)$ of the quasilinear
system
\beq\label{LWPPF}\begin{split}\pa_t\psi-\De\mu=f_1,\quad \mu=-\De\psi+\Phi'(\psi)-\la'(\psi)\vth,&\quad t\in J,\ x\in\Om,\\
\pa_t\left(b(\vth)+\la(\psi)\right)-\De\vth=f_2,&\quad t\in J,\ x\in\Om,\\
\pa_\nu\mu=g_1,\ \pa_\nu\psi=g_2,\ \pa_\nu\vth=g_3,&\quad t\in J,\ x\in\pa\Om,\\
\psi(0)=\psi_0,\ \vth(0)=\vth_0,&\quad t=0,\ x\in\Om.
\end{split}\eeq
To this end, we will apply Banach's fixed
point theorem. For this purpose let $p>(n+2)/2$, $p\ge 2$, $f_1,f_2\in X(T_0)$, $g_j\in
Y_j(0,T_0),\ j=1,2,3$, $\psi_0\in X_{\gamma}^1$ and $\vth_0\in X_\gamma^2$ be
given such that the compatibility conditions
$$\pa_\nu\De\psi_0-\pa_\nu\Phi'(\psi_0)+\pa_\nu(\la'(\psi_0)\vth_0)=-g_1|_{t=0},\
\pa_\nu\psi_0=g_2|_{t=0}\quad\text{and}\quad
\pa_\nu\vth_0=g_3|_{t=0}$$ are satisfied, whenever $p>5$, $p>5/3$ and $p>3$,
respectively. In the sequel we will assume that $\lambda,\phi\in C^{4-}(\R)$, $b\in C^{3-}(0,\infty)$ and $b'(s)>0$ for all $s>0$. Note that by the Sobolev embedding theorem we have $\vartheta_0\in C(\bar{\Omega})$ as well as $b'(\vartheta_0)\in C(\bar{\Omega})$. Since $\vartheta$ represents the inverse absolute temperature of the system, it is reasonable to assume $\vartheta_0(x)>0$ for all $x\in\bar{\Omega}$. Therefore, there exists a constant $\sigma>0$ such that $\vartheta_0(x),b'(\vartheta_0(x))\ge \sigma>0$ for all $x\in\bar{\Omega}$. We define $a_0(x):=1/b'(\vth_0(x))$, $\eta_1(x)=\la'(\psi_0(x))$ and $\eta_2(x)=a_0(x)\eta_1(x)$. By assumption, it holds that $a_0\in B_{pp}^{2-2/p}(\Omega)$, $\eta_1\in B_{pp}^{4-4/p}(\Omega)$ and $\eta_2\in B_{pp}^{2-2/p}(\Omega)$, cf. \cite[Section 4.6 \& Section 5.3.4]{RuSi96}.

Thanks to Theorem \ref{linthm} we
may define a pair of functions $(u^*,v^*)\in \E_1(T_0)$ as the solution
of the problem \beq\begin{split}\label{linfix}
\pa_tu^*+\De^2u^*+\De(\eta_1v^*)=f_1,&\quad t\in [0,T_0],\ x\in\Om,\\
\pa_tv^*-a_0\De v^*+\eta_2\pa_tu^*=a_0f_2,&\quad t\in [0,T_0],\ x\in\Om,\\
\pa_\nu\De u^*+\pa_\nu(\eta_1 v^*)=-g_1-e^{-B^2t}g_0,&\quad t\in [0,T_0],\ x\in\pa\Om,\\
\pa_\nu u^*=g_2,&\quad t\in [0,T_0],\ x\in\pa\Om,\\
\pa_\nu v^*=g_3,&\quad t\in [0,T_0],\ x\in\pa\Om,\\
u^*(0)=\psi_0,\ v^*(0)=\vth_0,&\quad t=0,\ x\in\Om,
\end{split}\eeq
where $B=-\De_{\partial\Omega}$ is the Laplace-Beltrami operator on $\partial\Omega$ and
$e^{-B^2t}$ is the analytic semigroup which is generated by $-B^2$.
Furthermore $g_0=0$ if $p<5$ and
$g_0=-g_1|_{t=0}-(\pa_\nu\De\psi_0+\pa_\nu(\eta_1\vth_0))$ if $p>5$.

Define a linear operator $\L:\!_0\E_1(T_0)\to\!_0\E_0(T_0)$ by
$$\L(u,v)=\vectfuenf{\pa_tu+\De^2u+\eta_1\De v}{\pa_tv-a_0\De v+\eta_2\pa_tu}{\pa_\nu\De u+\pa_\nu(\eta_1 v)}{\pa_\nu u}{\pa_\nu v}.$$
Then, by Theorem \ref{linthm}, the
operator $\L:\!_0\E_1(T_0)\to\!_0\E_0(T_0)$ is bounded and
bijective, hence an isomorphism with bounded inverse $\L^{-1}$.
For all $(u,v)\in\!_0\E_1(T)$ we set
$$G_1(u,v)=(\la'(\psi_0)-\la'(u)) v+\Phi'(u),$$
$$G_2(u,v)=(a_0\la'(\psi_0)-a(v)\la'(u))\pa_tu-(a_0-a(v))\De v-(a_0-a(v))f_2,$$
where $a(v(t,x))=1/b'(v(t,x))$ and $a_0=a(\vth_0)$. Lastly we
define a nonlinear mapping
$G:\E_1(T)\times\!_0\E_1(T)\to\!_0\E_0(T)$ by
$$G((u^*,v^*);(u,v))=\vectfuenf{\De G_1(u+u^*,v+v^*)}{G_2(u+u^*,v+v^*)}{\pa_\nu G_1(u+u^*,v+v^*)-\tilde{g}_0}{0}{0},$$
where $\tilde{g}_0=0$ if $p<5$ and $\tilde{g}_0=e^{-B^2 t}\pa_\nu
G_1(\psi_0,\vth_0)$ if $p>5$. Then it is easy to see that
$\psi=u+u^*\in E_1(T)$ and $\vth=v+v^*\in E_2(T)$ is a solution of \eqref{PF} if and
only if
$$\L(u,v)=G((u^*,v^*);(u,v))$$
or equivalently
$$(u,v)=\L^{-1}G((u^*,v^*);(u,v)).$$
In order to apply the contraction mapping principle we consider a
ball $\B_R=\B_R^1\times\B_R^2\subset\!_0\E_1(T)$, where
$R\in (0,1]$. Furthermore we define a mapping
$\calT:\B_R\to\!_0\E_1(T)$ by
$\calT(u,v)=\L^{-1}G((u^*,v^*);(u,v))$. We shall prove that
$\calT\B_R\subset\B_R$ and that $\calT$ defines a strict
contraction on $\B_R$. To this end we define the shifted ball
$\B_R(u^*,v^*)=\B_R^1(u^*)\times\B_R^2(v^*)\subset\E_1(T)$ by
$$\B_R(u^*,v^*)=\{(u,v)\in\E_1(T):(u,v)=(\tilde{u},\tilde{v})+(u^*,v^*),\ (\tilde{u},\tilde{v})\in\B_R\}.$$
To ensure that the mapping $G_2$ is well defined, we choose
$T_0>0$ and $R>0$ sufficiently small. This yields that all
functions $v\in\B_R^2(v^*)$ have only a small deviation from the
initial value $\vth_0$. To see this, write
$$|\vth_0(x)-v(t,x)|\le |\vth_0(x)-v^*(t,x)|+|v^*(t,x)-v(t,x)|\le
\mu(T)+R,$$ for all functions $v\in\B_R^2(v^*)$, where
$\mu=\mu(T)$ is defined by
$$\mu(T)=\max_{(t,x)\in[0,T]\times\Om}|v^*(t,x)-\vth_0(x)|.$$
Observe that $\mu(T)\to 0$ as $T\to 0$, by the continuity of $v^*$
and $\vth_0$. This in turn implies that $v(t,x)\ge\sigma/2>0$ and $b'(v(t,x))\ge\si/2>0$ for $(t,x)\in [0,T]\times\bar{\Omega}$ and all $v\in\B_R^2(v^*)$, with $T_0>0$, $R>0$ being sufficiently small.
Moreover, for all $v,\bar{v}\in\B_R^2(v^*)$ we obtain the
estimates \beq\label{lipscha1}|a(\vth_0(x))-a(v(t,x))|\le
C|\vth_0(x)-v(t,x)|\eeq and
\beq\label{lipscha2}|a(\bar{v}(t,x))-a(v(t,x))|\le
C|\bar{v}(t,x)-v(t,x)|,\eeq valid for all $(t,x)\in
[0,T]\times\bar{\Om}$, with some constant $C>0$, since $b'$ is locally
Lipschitz continuous.

The next proposition provides all the facts to show the
desired properties of the operator $\calT$.
\begin{prop}\label{fixpointest}
Let $n\in\N$ and $p>(n+2)/2$, $p\ge 2$, $b\in C^{2-}(0,\infty)$, $b'(s)>0$ for all $s>0$, $\lambda,\Phi\in C^{4-}(\R)$ and $\vth_0(x)>0$ for all $x\in\bar{\Omega}$. Then there
exists a constant $C>0$, independent of $T$, and functions
$\mu_j=\mu_j(T)$ with $\mu_j(T)\to 0$ as $T\to 0$, such that for
all $(u,v),(\bar{u},\bar{v})\in\B_R(u^*,v^*)$ the following
statements hold.
\begin{enumerate}
\setlength{\itemsep}{1ex}
\item $|\Delta G_1(u,v)-\Delta G_1(\bar{u},\bar{v})|_{X(T)}\le (\mu_1(T)+R)|(u,v)-(\bar{u},\bar{v})|_{\E_1(T)}$,
\item $|G_2(u,v)-G_2(\bar{u},\bar{v})|_{X(T)}\le C(\mu_2(T)+R)|(u,v)-(\bar{u},\bar{v})|_{\E_1(T)}$,
\item $|\pa_\nu G_1(u,v)-\pa_\nu G_1(\bar{u},\bar{v})|_{Y_1(T)}\le
C(\mu_3(T)+R)|(u,v)-(\bar{u},\bar{v})|_{\E_1(T)}$.
\end{enumerate}
\end{prop}
The proof is given in the Appendix.\vspace{0.2cm}

\noindent It is now easy to verify the self-mapping property of
$\calT$. Let $(u,v)\in \B_R$. By Proposition \ref{fixpointest}
there exists a function $\mu=\mu(T)$ with $\mu(T)\to 0$ as $T\to
0$ such that
\begin{align*}
|\calT(u,v)|_1&=|\L^{-1}G((u^*,v^*),(u,v))|_1\le |\L^{-1}||G((u^*,v^*),(u,v))|_0\\
&\le C(|G((u^*,v^*),(u,v))-G((u^*,v^*),(0,0))|_0+|G((u^*,v^*),(0,0))|_0)\\
&\le
C(|\Delta G_1(u+u^*,v+v^*)-\Delta G_1(u^*,v^*)|_{X(T)}\\
&\hspace{1cm}+|G_2(u+u^*,v+v^*)-G_2(u^*,v^*)|_{X(T)}\\
&\hspace{1cm}+|\pa_\nu G_1(u+u^*,v+v^*)-\pa_\nu G_1(u^*,v^*)|_{Y_1(T)}\\
&\hspace{1cm}+|G((u^*,v^*),(0,0))|_0)\\
&\le C(\mu(T)+R)|(u,v)|_1+|G((u^*,v^*),(0,0))|_0\\
&\le C(\mu(T)+R)R+|G((u^*,v^*),(0,0))|_0.
\end{align*}
Hence we see that $\calT\B_R\subset\B_R$ if $T$ and $R$ are
sufficiently small, since $G((u^*,v^*),(0,0))$ is a fixed
function. Furthermore for all $(u,v),(\bar{u},\bar{v})\in\B_R$ we
have
\begin{align*}
|\calT(u,v)-\calT(\bar{u},\bar{v})|_1&=|\L^{-1}(G((u^*,v^*),(u,v))-G((u^*,v^*),(\bar{u},\bar{v})))|_1\\
&\le |\L^{-1}||G((u^*,v^*),(u,v))-G((u^*,v^*),(\bar{u},\bar{v}))|_0\\
&\le
C(|\Delta G_1(u+u^*,v+v^*)-\Delta G_1(\bar{u}+u^*,\bar{v}+v^*)|_{X(T)}\\
&\hspace{1cm}+|\pa_\nu G_1(u+u^*,v+v^*)-\pa_\nu G_1(\bar{u}+u^*,\bar{v}+v^*)|_{Y_1(T)}\\
&\hspace{1cm}+|G_2(u+u^*,v+v^*)-G_2(\bar{u}+u^*,\bar{v}+v^*)|_{X(T)})\\
&\le C(\mu(T)+R)|(u,v)-(\bar{u},\bar{v})|_1.
\end{align*}
Thus $\calT$ is a strict contraction on $\B_R$, if $T$ and $R$ are
again small enough. Therefore we may apply the contraction mapping
principle to obtain a unique fixed point
$(\tilde{u},\tilde{v})\in\B_R$ of $\calT$. In other words the
pair $(\psi,\vth)=(\tilde{u}+u^*,\tilde{v}+v^*)\in\E_1(T)$ is
the unique local solution of \eqref{PF}. We summarize the
preceding calculations in
\begin{thm}\label{LWPsolPF}
Let $n\in\N$, $p>(n+2)/2$, $p\ge 2$, $p\neq 3,5$, $b\in C^{3-}(0,\infty)$, $b'(s)>0$ for all $s>0$ and let $\lambda,\Phi\in C^{4-}(\R)$. Then there exists an interval
$J=[0,T]\subset [0,T_0]=J_0$ and a unique solution $(\psi,\vth)$ of
\eqref{PF} on $J$, with
$$\psi\in H_p^1(J;L_p(\Om))\cap L_p(J;H_p^4(\Om))$$
and
$$\vth\in H_p^1(J;L_p(\Om))\cap L_p(J;H_p^2(\Om)),\quad \vth(t,x)>0\ \mbox{for all}\ (t,x)\in J\times\bar{\Omega},$$
provided the data are subject to the following conditions.
\begin{enumerate}
\item $f_1,f_2\in L_p(J_0\times\Om)$,
\item $g_1\in W_p^{1/4-1/4p}(J_0;L_p(\partial\Omega))\cap
L_p(J_0;W_p^{1-1/p}(\partial\Omega))$,
\item $g_2\in W_p^{3/4-1/4p}(J_0;L_p(\partial\Omega))\cap
L_p(J_0;W_p^{3-1/p}(\partial\Omega))$,
\item $g_3\in W_p^{1/2-1/2p}(J_0;L_p(\partial\Omega))\cap
L_p(J_0;W_p^{1-1/p}(\partial\Omega))$,
\item $\psi_0\in B_{pp}^{4-4/p}(\Om)$, $\vth_0\in B_{pp}^{2-2/p}(\Om)$,
\item $\pa_\nu\De \psi_0-\pa_\nu\Phi'(\psi_0)+\pa_\nu(\la'(\psi_0) \vth_0)=-g_1|_{t=0},$ if
$p>5$,
\item $\pa_\nu\psi_0=g_2|_{t=0}$, $\pa_\nu \vth_0=g_3|_{t=0},$ if $p>3$,
\item $\vth_0(x)>0$ for all $x\in\bar{\Omega}$.
\end{enumerate}
The solution depends continuously on the given data and if the
data are independent of $t$, the map $(\psi_0,\vth_0)\mapsto
(\psi,\vth)$ defines a local semiflow on the natural (nonlinear) phase
manifold
$$\calM_p:=\{(\psi_0,\vth_0)\in B_{pp}^{4-4/p}(\Om)\times B_{pp}^{2-2/p}(\Om): \psi_0\ \text{and}\ \vth_0\ \text{satisfy}\ 6.-8.\}.$$
\end{thm}

\section{Global Well-Posedness}\label{SECPFGWP}
In this section we will investigate the global existence of the
solution to the conserved Penrose-Fife type system
\beq\label{PFGWP}\begin{split}\pa_t\psi-\De\mu=0,\quad \mu=-\De\psi+\Phi'(\psi)-\la'(\psi)\vth,&\quad t>0,\ x\in\Om,\\
\pa_t\left(b(\vth)+\la(\psi)\right)-\De\vth=0,&\quad t>0,\ x\in\Om,\\
\pa_\nu\mu=0,\ \pa_\nu\psi=0,\ \pa_\nu\vth=0,&\quad t>0,\ x\in\pa\Om,\\
\psi(0)=\psi_0,\ \vth(0)=\vth_0,&\quad t=0,\ x\in\Om,
\end{split}\eeq
with respect to time if the spatial dimension $n$ is less or equal to 3.
Note that the boundary conditions are equivalent to
$\pa_\nu\vth=\pa_\nu\psi=\pa_\nu\De\psi=0$. A successive application of Theorem \ref{LWPsolPF} yields
a maximal interval of existence $J_{\max}=[0,T_{\max})$ for the
solution $(\psi,\vth)\in E_1(T)\times E_2(T)$ of \eqref{PFGWP}, where $T\in(0,T_{\max})$. In the
sequel we will make use of the following assumptions.
\begin{itemize}
\setlength{\itemsep}{0.5ex}
\item[\textbf{(H1)}] $\Phi\in C^{4-}(\mathbb{R})$ and there exist some constants $c_j>0$, $\ga\ge 0$
such that
$$\Phi(s)\ge-\frac{\eta}{2}s^2-c_1,\ |\Phi'''(s)|\le c_2(1+|s|^{\ga}),$$
for all $s\in\R$, where $\eta<\la_1$ with $\la_1$ being the
smallest nontrivial eigenvalue of the negative Laplacian on $\Om$
with Neumann boundary conditions and $\ga<3$ if $n=3$.
\item[\textbf{(H2)}] $\lambda\in C^{4-}(\mathbb{R})$ and $\la'',\la'''\in L_\infty(\R)$. In particular, there is a constant $c>0$ such that
$|\la'(s)|\le c(1+|s|)$ for all $s\in\R$.
\item[\textbf{(H3)}] $b\in C^{3-}((0,\infty))$, $b'(s)>0$ on $(0,\infty)$ and there is
a constant $\ka>1$ such that $$\frac{1}{\kappa}\le\vartheta(t,x)\le\kappa$$ on
$J_{\max}\times\Om$. In particular, there exists $\sigma>1$ such that
    $$\frac{1}{\sigma}\le b'(\vartheta(t,x))\le\sigma,$$
on $J_{\max}\times\Omega$.
\end{itemize}
\emph{Remark:} Condition (H1) is certainly fulfilled, if $\Phi$ is
a polynomial of degree $2m$, $m<3$.\vspace{0.25cm}

\noindent We prove global well-posedness with respect to time by contradiction. For
this purpose, assume that $T_{\max}<\infty$. Multiply
$\pa_t\psi=\De\mu$ by $\mu$ and integrate by parts to the result
\beq\label{eqpsi}\difft\left(\frac{1}{2}|\nabla\psi|_{2}^2+\int_\Om\Phi(\psi)\
dx\right)+|\nabla\mu|_{2}^2-\int_\Om\la'(\psi)\vth\pa_t\psi\
dx=0.\eeq Next we multiply $\eqref{PFGWP}_2$ by $\vth$ and
integrate by parts. This yields \beq\label{eqtheta}\int_\Om\vth
b'(\vth)\pa_t\vth\
dx+|\nabla\vth|_{2}^2+\int_\Om\la'(\psi)\vth\pa_t\psi\ dx=0.\eeq
Set $\beta'(s)=sb'(s)$ and add \eqref{eqpsi} to \eqref{eqtheta} to
obtain the equation
\begin{align}\label{eneq}\begin{split}\difft\Big(\frac{1}{2}|\nabla\psi|_{2}^2+
\int_\Om\Phi(\psi)\ dx+\int_\Om\beta(\vth)\
dx\Big)+|\nabla\mu|_{2}^2+|\nabla\vth|_{2}^2=0.\end{split}\end{align}
Integrating \eqref{eneq} with respect to $t$, we obtain
\beq\label{PFenineq2}
E(\psi(t),\vth(t))+\int_0^t\left(|\nabla\mu(s)|_{2}^2+|\nabla\vth(s)|_{2}^2\right)\ dt=E(\psi_0,\vth_0),\eeq
for all $t\in J_{\max}$, where
$$E(u,v):=\frac{1}{2}|\nabla u|_{2}^2+\int_\Om\Phi(u)\
dx+\int_\Om\beta(v)\ dx.$$ It follows from (H1) and the
Poincar\'{e}-Wirtinger inequality that
\begin{multline*}\frac{\ep}{2}\int_\Om|\nabla\psi(t)|^2\
dx+\frac{1-\ep}{2}\int_\Om|\nabla\psi(t)|^2\
dx+\int_\Om\Phi(\psi(t))\ dx\\
\ge\frac{\ep}{2}\int_\Om|\nabla\psi(t)|^2\ dx+
\frac{(1-\ep)\la_1-\eta}{2}|\psi(t)|_2^2-c_1|\Om|-\frac{\la_1}{2|\Om|}\left(\int_{\Om}\psi_0\
dx\right),\end{multline*} since by equation $\pa_t\psi=\De\mu$ and
the boundary condition $\pa_\nu\mu=0$, it holds that
$$\int_\Om\psi(t,x)\ dx\equiv \int_\Om\psi_0(x)\ dx,\quad t\in J_{\max}.$$
Hence for a sufficiently small $\ep>0$ we obtain the a priori
estimates \beq\label{enest2}\psi\in
L_\infty(J_{\max};H_2^1(\Om))\quad\text{and}\quad
|\nabla\mu|,|\nabla\vth|\in L_2(J_{\max};L_2(\Om)),\eeq since
$\beta(\vth(t,x))$ is uniformly bounded on $J_{\max}\times\Om$, by
(H3). However, things are more involved for higher order estimates.
Here we have the following result.
\begin{prop}\label{estpsit}
Let $n\le 3$, $p>(n+2)/2$, $p\ge 2$ and let $(\psi,\vth)$ be the
maximal solution of \eqref{PFGWP} with initial value
$\psi_0\in B_{pp}^{4-4/p}(\Omega)$ and $\vth_0\in B_{pp}^{2-2/p}(\Omega)$. Suppose furthermore $b\in C^{3-}(0,\infty)$, $b'(s)>0$ for all $s>0$, $\lambda,\Phi\in C^{4-}(\R)$ and let (H1)-(H3) hold.

Then $\psi\in L_\infty(J_{\max}\times\Om)$ and $\vth\in
H_2^1(J_{\max};L_2(\Om))\cap L_\infty(J_{\max};H_2^1(\Om))$.
Moreover, it holds that $\pa_t\psi\in L_r(J_{\max}\times\Om)$,
where $r:=\min\{p,2(n+4)/n\}$.
\end{prop}
\bpr The proof is given in the Appendix.

\epr \noindent Define the new function $u=b(\vth)$. Then u
satisfies the nonautonomous linear differential equation in divergence form
\beq\label{eqhoelder}\pa_t u-\diver (a(t,x)\nabla u)=f,\eeq
subject to the boundary and initial conditions $\pa_\nu u=0$ and
$u(0)=b(\vth_0)=:u_0$, where $a(t,x):=1/b'(\vth(t,x))$ and
$f:=-\la'(\psi)\pa_t\psi$. With (H3), the regularity of $\vth$ from
Proposition \ref{estpsit} carries over to the function $u$; in particular $u_0\in B_{pp}^{2-2/p}(\Omega)$. This yields, that $u$ is a
\emph{weak solution} of \eqref{eqhoelder} in the sense of
\textsc{Lieberman} \cite{Lieb96} \& \textsc{DiBenedetto} \cite{DiBen93}, and $u$ is bounded by (H3).

Furthermore, by (H3)
$$0<\frac{1}{\sigma}\le a(t,x)\le \sigma<\infty,$$
for all $(t,x)\in J_{\max}\times\Om$. Note that by Proposition \ref{estpsit} it holds that $f=-\la'(\psi)\pa_t\psi\in
L_r(J_{\max}\times\Om)$, $r:=\min\{p,2(n+4)/n\}$. Consider the case
$r=2(n+4)/n$. Then it can be readily checked that
$$\frac{n+2}{2}<\frac{2(n+4)}{n}=r$$
provided $n\le 5$. It follows from \textsc{Lieberman} \cite{Lieb96} \& \textsc{DiBenedetto} \cite{DiBen93} that there exists a real number $\al\in
(0,1/2)$ such that $u\in C^{\al,2\al}(\overline{\Omega_{T_{\max}}})$,
provided $f\in L_p(J_{\max}\times\Om)$ and $p>(n+2)/2$. Here
$C^{\al,2\al}(\overline{\Omega_{T_{\max}}})$ is defined as
$$C^{\al,2\al}(\overline{\Om_{T_{\max}}}):=\{v\in
C(\overline{\Om_{T_{\max}}}):\sup_{(t,x),(s,y)\in
\overline{\Om_{T_{\max}}}}\frac{|v(t,x)-v(s,y)|}{|t-s|^\al+|x-y|^{2\al}}<\infty\}.$$
and we have set $\Om_{T_{\max}}=(0,T_{\max})\times\Om$. The properties of the function $b$
then yield that $\vth=b^{-1}(u)\in C^{\alpha,2\alpha}(\overline{\Omega_{T_{\max}}})$.
In a next step we solve the
initial-boundary value problem
\begin{align}\label{PFGWPEQVTH}\begin{split}
\pa_t\vth-a(t,x)\De\vth=g,&\quad
t\in J_{\max},\ x\in\Om,\\
\pa_\nu\vth=0,&\quad t\in J_{\max},\ x\in\pa\Om,\\
\vth(0)=\vth_0,&\quad t=0,\ x\in\Om, \end{split}\end{align} with
$g:=-a(t,x)\la'(\psi)\pa_t\psi\in L_r(J_{\max}\times\Om)$ and
$r=2(n+4)/n>(n+2)/2$. By \cite[Theorem 2.1]{DHP07} we obtain
$$\vth\in H_r^1(J_{\max};L_r(\Om))\cap L_r(J_{\max};H_r^2(\Om)),$$
of \eqref{PFGWPEQVTH}, since
$$\vth_0\in
B_{pp}^{2-2/p}(\Om)\hookrightarrow B_{rr}^{2-2/r}(\Om),\quad p\ge
r.$$ At this point we use equation \eqref{psit} from the proof of
Proposition \ref{estpsit} to conclude $\pa_t\psi\in
L_s(J_{\max}\times\Om)$, with $s=\min\{p,q\}$ where $q$ is
restricted by
$$\frac{1}{q}\ge\frac{1}{r}-\frac{2}{n+4}.$$
For the case $r=2(n+4)/n$, this yields
$$\frac{1}{q}\ge\frac{n-4}{2(n+4)},$$
i.e. $q$ may be arbitrarily large in case $n\le 3$ and we may set
$s=p$. Now we solve \eqref{PFGWPEQVTH} again, this time with $g\in
L_p(J_{\max}\times\Om)$, to obtain
$$\vth\in H_p^1(J_{\max};L_p(\Om))\cap L_p(J_{\max};H_p^2(\Om))$$
and therefore $\vth(T_{\max})\in B_{pp}^{2-2/p}(\Om)$ is well
defined. Next, consider the equation
$$\pa_t\psi+\De^2\psi=\De\Phi'(\psi)-\De(\la'(\psi)\vth),$$
subject to the initial and boundary conditions $\psi(0)=\psi_0$
and $\pa_\nu\psi=\pa_\nu\De\psi=0$. By maximal $L_p$-regularity there
exists a constant $M=M(J_{\max})>0$ such that
    \beq\label{MRestPF}|\psi|_{E_1(T)}\le
    M(1+|\De\Phi'(\psi)|_{X(T)}+|\De(\la'(\psi)\vth)|_{X(T)}).
    \eeq
for each $T\in J_{\max}$. Since $\vth\in E_2(T_{\max})$ we may apply \cite[Lemma 4.1]{PrWi06} to the result
    \begin{equation}\label{MRestPF2}
    |\De\Phi'(\psi)|_{X(T)}+|\De(\la'(\psi)\vth)|_{X(T)}\le C(1+|\psi|_{E_1(T)}^\delta),
    \end{equation}
with some $\delta\in (0,1)$ and $C>0$ being independent of $T\in J_{\max}$. Combining \eqref{MRestPF} with \eqref{MRestPF2}, we obtain the estimate
    $$|\psi|_{E_1(T)}\le C(1+|\psi|_{E_1(T)}^\delta),$$
which in turn yields that $|\psi|_{E_1(T)}$ is bounded as $T\nearrow T_{\max}$, since $\delta\in (0,1)$. Therefore the value $\psi(T_{\max})\in B_{pp}^{4-4/p}(\Omega)$ is well defined and we may continue the solution $(\psi,\vth)$ beyond the point $T_{\max}$, contradicting the assumption that $J_{\max}=[0,T_{\max})$ is the maximal interval of existence.
We summarize these considerations in
\begin{thm}\label{GWPsolPF}
Let $n\le 3$, $p>(n+2)/2$, $p\ge 2$ and $p\neq 3,5$. Assume that (H1)-(H3) hold. Then for each $T_0>0$ there exists a unique
solution
$$\psi\in H_p^1(J_0;L_p(\Om))\cap L_p(J_0;H_p^4(\Om))=E_1(T_0)$$
and
$$\vth\in H_p^1(J_0;L_p(\Om))\cap L_p(J_0;H_p^2(\Om))=E_2(T_0),$$
of \eqref{PF}, provided the data are subject to the following conditions.
\be
\item $\psi_0\in B_{pp}^{4-4/p}(\Om)$, $\vth_0\in B_{pp}^{2-2/p}(\Om)$;
\item $\pa_\nu\De \psi_0=0,$ if
$p>5$, $\pa_\nu\psi_0=0$;
\item $\pa_\nu \vth_0=0,$ if $p>3$, $\vth_0(x)>0$ for all $x\in\bar{\Omega}$.
\ee The solution depends continuously on the given data and the
map $(\psi_0,\vth_0)\mapsto (\psi,\vth)$ defines a semiflow
on the natural phase manifold
$$\calM_p:=\{(\psi_0,\vth_0)\in B_{pp}^{4-4/p}(\Om)\times B_{pp}^{2-2/p}(\Om): \psi_0\ \text{and}\ \vth_0\ \text{satisfy}\ 2.\ \&\ 3.\}.$$
\end{thm}

\section{Asymptotic Behavior}\label{asymbeh}

Let $n\le3$. In the following we will investigate the asymptotic
behavior of global solutions of the homogeneous system
\beq\label{PFasym}\begin{split}\pa_t\psi-\De\mu=0,\quad \mu=-\De\psi+\Phi'(\psi)-\la'(\psi)\vth,&\quad t>0,\ x\in\Om,\\
\pa_t\left(b(\vth)+\la(\psi)\right)-\De\vth=0,&\quad t>0,\ x\in\Om,\\
\pa_\nu\mu=0,&\quad t>0,\ x\in\pa\Om,\\
\pa_\nu\psi=0,&\quad t>0,\ x\in\pa\Om,\\
\pa_\nu\vth=0,&\quad t>0,\ x\in\pa\Om,\\
\psi(0)=\psi_0,\ \vth(0)=\vth_0,&\quad t=0,\ x\in\Om,
\end{split}\eeq
as $t\to \infty$. To this end let $(\psi_0,\vth_0)\in\calM_p$, $p>(n+2)/2$, $p\ge 2$ and denote
by $(\psi(t),\vth(t))$ the unique global solution of
\eqref{PFasym}. In the sequel we will make use of the following
assumptions.
\begin{itemize}
\setlength{\itemsep}{0.5ex}
\item[\textbf{(H4)}]$b\in C^{3-}((0,\infty))$, $b'(s)>0$ on $(0,\infty)$ and there is
a constant $\ka>1$ such that $$\frac{1}{\kappa}\le\vartheta(t,x)\le\kappa$$ on
$J_{\max}\times\Om$. In particular, there exists $\sigma>1$ such that
    $$\frac{1}{\sigma}\le b'(\vartheta(t,x))\le\sigma,$$
on $J_{\max}\times\Omega$.
\item[\textbf{(H5)}] The functions $\Phi$, $\la$ and $b$ are real
analytic on $\R$.
\end{itemize}
We remark that assumption (H4) is identical to (H3) for a global solution. We stated it here for the sake of readability.

Note that the boundary conditions $\eqref{PFasym}_{3,5}$ yield
$$\int_\Om \psi(t,x)\ dx\equiv\int_\Om \psi_0(x)\
dx,$$ and $$\int_\Om(b(\vth(t,x))+\la(\psi(t,x)))\ dx\equiv
\int_\Om(b(\vth_0(x))+\la(\psi_0(x)))\ dx.$$ Replacing $\psi$ by
$\tilde{\psi}=\psi-c$, where $c:=\frac{1}{|\Om|}\int_\Om
\psi_0(x)\ dx$ we see that $\int_\Om\tilde{\psi}\ dx \equiv 0$, if
$\Phi(s)$ and $\la(s)$ are replaced by $\tilde{\Phi}(s)=\Phi(s+c)$
and $\tilde{\la}(s)=\la(s+c)$, respectively. Similarly we can
achieve that
$$\int_\Om(b(\vth(t,x))+\la(\psi(t,x)))\ dx\equiv 0,$$
by a shift of $\la$, to be precise $\bar{\la}(s):=\la(s)-d$,
where
$$d:=\frac{1}{|\Om|}\int_\Om(b(\vth_0(x))+\la(\psi_0(x)))\ dx.$$
With these modifications of the data we obtain the constraints
\beq\label{sidecond}\int_\Om \psi(t,x)\
dx\equiv0\quad\text{and}\quad
\int_\Om(b(\vth(t,x))+\la(\psi(t,x)))\ dx\equiv 0.\eeq Recall from
Section \ref{SECPFGWP} the energy functional
$$E(u,v)=\frac{1}{2}|\nabla u|_{2}^2+\int_\Om\Phi(u)\
dx+\int_\Om\beta(v)\ dx,$$ defined on the energy space
$V=V_1\times V_2$, where
$$V_1:=\left\{u\in H_2^1(\Om):\int_\Om u\
dx=0\right\},\qquad V_2:=H_2^r(\Om),\ r\in (n/4,1).$$ and $V$ is
equipped with the canonical norm
$|(u,v)|_V:=|u|_{H_2^1(\Om)}+|v|_{H_2^r( \Om)}$. It is convenient
to embed $V$ into a Hilbert space $H=H_1\times H_2$ where
$$H_1:=\left\{u\in L_2(\Om):\int_\Om u\ dx=0\right\}\quad\text{and}\quad H_2:=L_2(\Om).$$
\begin{prop}\label{relcomp}
Let $(\psi,\vth)\in E_1\times E_2$ be a global solution of \eqref{PFasym} and
assume (H1)-(H4). Then
\begin{enumerate}
\item $\psi\in L_\infty(\R_+;H_p^{2s}(\Om)),\ s\in [0,1),\ p\in (1,\infty),\ \pa_t\psi\in
L_2(\R_+\times\Om)$;
\item $\vth\in L_\infty(\R_+;H_2^1(\Om)),\ \pa_t\vth\in
L_2(\R_+\times\Om)$.
\end{enumerate}
In particular the orbits $\psi(\R_+)$ and $\vth(\R_+)$ are
relatively compact in $H_2^1(\Om)$ and $H_2^{r}(\Om)$,
respectively, where $r\in [0,1)$.\end{prop} \bpr Assertions 1 \& 2
follow directly from (H1)-(H4) and the proof of Proposition
\ref{estpsit}, which is given in the Appendix. Indeed, one may
replace the interval $J_{\max}$ by $\R_+$, since the operator
$-A^2=-\De_N^2$ generates an exponentially stable, analytic semigroup
$e^{-A^2 t}$ in the space
$$\X_p:=\{u\in L_p(\Om):\int_\Om u\ dx=0\}$$
with domain
$$D(A^2)=\{u\in H_p^4(\Om)\cap\X_p:\pa_\nu u=\pa_\nu\De u=0\ \text{on}\
\partial\Omega\}.$$
\epr

By Assumption (H4), there exists some bounded interval
$J_\vth\subset\R_+$ with $\vth(t,x)\in J_\vth$ for all $t\ge 0,\
x\in\Om$. Therefore we may modify the nonlinearities $b$ and
$\beta$ outside $J_\vth$ in such a way that $b,\beta\in
C_b^{3-}(\R)$.

Unfortunately the energy functional $E$ is not yet the right one for our purpose,
since we have to include the nonlinear constraint
$$\int_\Om(\la(\psi)+b(\vth))\ dx=0,$$
into our considerations. The linear constraint $\int_\Om\psi\
dx=0$ is part of the definition of the space $H_1$. For the
nonlinear constraint we use a functional of Lagrangian type
which is given by
$$L(u,v)=E(u,v)-\overline{v} F(u,v),$$
defined on $V$, where $F(u,v):=\int_\Om(\la(u)+b(v))\ dx$ and $\bar{w}=\frac{1}{|\Om|}\int_\Om w\ dx$ for a
function $w\in L_1(\Om)$. Concerning the differentiability of $L$ we
have the following result.
\begin{prop}\label{diffL}
Under the conditions (H1)-(H4), the functional $L$ is twice continuously Fr\'{e}chet
differentiable on $V$ and the derivatives are given by
\begin{multline}\label{1stder}
\langle L'(u,v),(h,k)\rangle_{V^*,V}=\\
\langle E'(u,v),(h,k)\rangle_{V^*,V}-\overline{k}F(u,v)-\overline{v}\langle
F'(u,v),(h,k)\rangle_{V^*,V}
\end{multline}
and
\begin{multline}\label{2ndder}
\langle L''(u,v)(h_1,k_1),(h_2,k_2)\rangle_{V^*,V}=
\langle E''(u,v)(h_1,k_1),(h_2,k_2)\rangle_{V^*,V}-\\
\overline{k_1}\langle F'(u,v),(h_2,k_2)\rangle_{V^*,V}-\overline{k_2}\langle
F'(u,v),(h_1,k_1)\rangle_{V^*,V}-\\
\overline{v}\langle
F''(u,v)(h_1,k_1),(h_2,k_2)\rangle_{V^*,V},
\end{multline}
where $(h,k),(h_j,k_j)\in V,\ j=1,2$, and
\begin{align*}
\langle E'(u,v),(h,k)\rangle_{V^*,V}=\int_\Om\nabla u\nabla h\
dx+\int_\Om\Phi'(u) h\ dx+\int_\Om\beta'(v) k\ dx,
\end{align*}
\begin{multline*}
\langle E''(u,v)(h_1,k_1),(h_2,k_2)\rangle_{V^*,V}=\\
\int_\Om\nabla
h_1\nabla h_2\ dx+\int_\Om\Phi''(u) h_1h_2\ dx+\int_\Om\beta''(v)
k_1k_2\ dx,
\end{multline*}
\begin{align*}
\langle F'(u,v),(h,k)\rangle_{V^*,V}=\int_\Om\la'(u) h\
dx+\int_\Om b'(v) k\ dx
\end{align*}
and
\begin{align*}
\langle
F''(u,v)(h_1,k_1),(h_2,k_2)\rangle_{V^*,V}=\int_\Om\la''(u)
h_1h_2\ dx+\int_\Om b''(v) k_1k_2\ dx.
\end{align*}
\end{prop}

\bpr We only consider the first derivative, the second one is
treated in a similar way. Since the bilinear form
\beq\label{bilinform}a(u,v):=\int_\Om\nabla u(x)\nabla v(x)\
dx\eeq defined on $V_1\times V_1$ is bounded and symmetric, the
first term in $E$ is twice continuously Fr\'{e}chet
differentiable. For the functional
$$G_1(u):=\int_\Om\Phi(u)\ dx,\quad u\in V_1,$$
we argue as follows. With $u,h\in V_1$ it holds that
\begin{align*}
\Phi(u(x)+h(x))-\Phi(u(x))&-\Phi'(u(x))h(x)\\
&=\int_0^1\frac{d}{dt}\ \Phi(u(x)+th(x))\ dt-\int_0^1\Phi'(u(x))h(x)\ dt\\
&=\int_0^1\Big(\Phi'(u(x)+th(x))-\Phi'(u(x))\Big)h(x)\ dt\\
&=\int_0^1\int_0^t\frac{d}{ds}\ \Phi'(u(x)+sh(x))h(x)\ ds\ dt\\
&=\int_0^1\int_0^t\Phi''(u(x)+sh(x))h(x)^2\ ds\ dt\\
&=\int_0^1\Phi''(u(x)+sh(x))h(x)^2(1-s)\ ds.
\end{align*}
From the growth condition (H1), H\"{o}lder's inequality and the
Sobolev embedding theorem it follows that
\begin{align*}
\Big|\int_\Om\Big(\Phi(u(x)+h(x))&-\Phi(u(x))-\Phi'(u(x))h(x)\Big)\ dx\Big|\\
&\le C\int_\Om(1+|u(x)|^4+|h(x)|^4)|h(x)|^2\ dx\\
&\le C(1+|u|_{6}^4+|h|_{6}^4)|h|_{6}^2\\
&\le C(1+|u|_{V_1}^4+|h|_{V_1}^4)|h|_{V_1}^2.
\end{align*}
This proves that $G_1$ is Fr\'{e}chet differentiable and also
$G_1'(u)=\Phi'(u)\in L_{6/5}(\Om)\hookrightarrow V_1^*$. The next
step is the proof of the continuity of $G_1':V_1\to V_1^*$. We
make again use of (H1), the H\"{o}lder inequality and the Sobolev
embedding theorem to obtain
\begin{align*}|G_1'(u)&-G_1'(\bar{u})|_{V_1^*}\\
&\le
C\left(\int_\Om|\Phi'(u(x))-\Phi'(\bar{u}(x))|^{\frac{6}{5}}\
dx\right)^\frac{5}{6}\\
&\le
C\left(\int_\Om\int_0^1|\Phi''(tu(x)+(1-t)\bar{u}(x))|^\frac{6}{5}|u(x)-\bar{u}(x)|^\frac{6}{5}\
dt\ dx\right)^\frac{5}{6}\\
&\le C
\left(\int_\Om(1+|u(x)|^\frac{24}{5}+|\bar{u}(x)|^\frac{24}{5})|u(x)-\bar{u}(x)|^\frac{6}{5}\
dx\right)^\frac{5}{6}\\
&\le C\left(\int_\Om(1+|u(x)|^6+|\bar{u}(x)|^6)\
dx\right)^\frac{2}{3}\left(\int_\Om|u(x)-\bar{u}(x)|^{6}\right)^{\frac{1}{6}}\\
&\le C(1+|u|_{V_1}^4+|\bar{u}|_{V_1}^4)|u-\bar{u}|_{V_1}.
\end{align*}
Actually this proves that $G_1'$ is even locally Lipschitz
continuous on $V_1$. The Fr\'{e}chet differentiability of $G_1'$
and the continuity of $G_1''$ can be proved in an analogue way.
The fundamental theorem of differential calculus and the Sobolev
embedding theorem yield the estimate
\begin{multline*}
|\Phi'(u+h)-\Phi'(u)-\Phi''(u)h|_{V_1^*}\\
\le C\left(\int_\Om\int_0^1|\Phi'''(u(x)+sh(x))|^{\frac{6}{5}}|h(x)|^{\frac{12}{5}}\
ds\ dx\right)^{\frac{5}{6}}.
\end{multline*}
We apply Assumption (H1) and H\"{o}lder's inequality to the result
\begin{align*}
|\Phi'(u+h)-\Phi'(u)&-\Phi''(u)h|_{V_1^*}\\
&\le C\left(\int_\Om(1+|u(x)|^\frac{18}{5}+|h(x)|^\frac{18}{5})|h(x)|^\frac{12}{5}\
dx\right)^{\frac{5}{6}}\\
&\le C\left(\int_\Om(1+|u(x)|^6+|h(x)|^6)\
dx\right)^{\frac{1}{2}}\left(\int_\Om|h(x)|^6\
dx\right)^{\frac{1}{3}}\\
&=C(1+|u|_{V_1}^3+|h|_{V_1}^3)|h|_{V_1}^2.
\end{align*}
Hence the Fr\'{e}chet derivative is given by the multiplication
operator $G_1''(u)$ defined by $G_1''(u)v=\Phi''(u)v$ for all
$v\in V_1$ and $\Phi''(u)\in L_{3/2}(\Om)$. We will omit the proof
of continuity of $G_1''$. The way to show the $C^2$-property of
the functional
$$G_2(u):=\int_\Om\la(u(x))\ dx,\quad u\in V_1,$$
is identical to the one above, by Assumption (H2). Concerning the
$C^2$-differentiability of the functionals
$$G_3(v):=\int_\Om\beta(v(x))\ dx\quad\text{and}\quad
G_4(v):=\int_\Om b(v(x))\ dx,\quad v\in V_2,$$ one may adopt the
proof for $G_1$ and $G_2$. In fact, this time it is easier, since
$\beta$ and $b$ are assumed to be elements of the space
$C_b^{3-}(\R)$, however one needs the assumption $r\in (n/4,1)$. We will skip the details.
Finally the product rule of differentiation yields that $L$ is
twice continuously Fr\'{e}chet differentiable on $V_1\times V_2$.

\epr

The corresponding stationary system to
\eqref{PFasym} will be of importance for the forthcoming
calculations. Setting all time-derivatives in \eqref{PFasym} equal
to 0 yields
$$\De\mu=0\quad\text{and}\quad\De\vth=0,$$
subject to the boundary conditions $\pa_\nu\mu=\pa_\nu\vth=0$.
Thus we have $\mu\equiv\mu_\infty=const$,
$\vth\equiv\vth_\infty=const$ and there remains the nonlinear
elliptic problem of second order \beq\label{statsys}
\begin{cases}
-\De\psi_\infty+\Phi'(\psi_\infty)-\la'(\psi_\infty)\vth_\infty=\mu_\infty,\quad
x\in\Om,\\
\pa_\nu\psi_\infty=0,\quad x\in\pa\Om,
\end{cases}
\eeq with the constraints \eqref{sidecond} for the unknowns
$\psi_\infty$ and $\vth_\infty$. The following proposition
collects some properties of the functional $L$ and the $\om$-limit
set
    \begin{multline*}
    \om(\psi,\vth):=\{(\ph,\th)\in V_1\times V_2:\exists\ (t_n)\nearrow \infty\ \mbox{s.t.}\\
    (\psi(t_n),\vth(t_n))\to (\ph,\th)\ \mbox{in}\ V_1\times V_2\}.
    \end{multline*}
\begin{prop}\label{propomlim}
Under Hypotheses (H1)-(H4) the following assertions are true.
\begin{enumerate}
\item The $\om$-limit set is nonempty, connected and compact.
\item Each point $(\psi_\infty,\vth_\infty)\in\om(\psi,\vth)$ is a strong solution of the
stationary problem \eqref{statsys}, where
$\vth_\infty,\mu_\infty=const$ and $(\psi_\infty,\vth_\infty)$
satisfies the constraints \eqref{sidecond} for the unknowns
$\vth_\infty,\mu_\infty$.
\item The functional $L$ is constant on $\om(\psi,\vth)$ and each
point $(\psi_\infty,\vth_\infty)\in\om(\psi,\vth)$ is a critical
point of $L$, i.e. $L'(\psi_\infty,\vth_\infty)=0$ in $V^*$.
\end{enumerate}
\end{prop}
\bpr The fact that $\om(\psi,\vth)$ is nonempty, connected and
compact follows from Proposition \ref{relcomp} and some well-known
facts in the theory of dynamical systems.

Now we turn to 2. Let $(\psi_\infty,\vth_\infty)\in
\om(\psi,\vth)$. Then there exists a sequence
$(t_n)\nearrow+\infty$ such that
$(\psi(t_n),\vth(t_n))\to(\psi_\infty,\vth_\infty)$ in $V$ as
$n\to\infty$. Since $\pa_t\psi,\pa_t\vth\in L_2(\R_+\times\Om)$ it
follows that $\psi(t_n+s)\to\psi_\infty$ and
$\vth(t_n+s)\to\vth_\infty$ in $L_2(\Om)$ for all $s\in [0,1]$ and
by relative compactness also in $V$. This can be seen as follows.
\begin{align*}
|\psi(t_n+s)-\psi_\infty|_{2}&\le
|\psi(t_n+s)-\psi(t_n)|_{2}+|\psi(t_n)-\psi_\infty|_{2}\\
&\le \int_{t_n}^{t_n+s}|\pa_t\psi(t)|_{2}\
dt+|\psi(t_n)-\psi_\infty|_{2}\\
&\le s^{1/2}\left(\int_{t_n}^{t_n+s}|\pa_t\psi(t)|_{2}^2\
dt\right)^{1/2}+|\psi(t_n)-\psi_\infty|_{2}.
\end{align*}
Then, for $t_n\to\infty$ this yields $\psi(t_n+s)\to\psi_\infty$ for
all $s\in[0,1]$. The proof for $\vth$ is the same. Integrating
\eqref{eneq} with $f_1=f_2=0$ from $t_n$ to $t_n+1$ we obtain
    \begin{multline*}
    E(\psi(t_n+1),\vth(t_n+1))-E(\psi(t_n),\vth(t_n))\\
    +\int_0^1\int_\Om\left(|\nabla\mu(t_n+s,x)|^2+|\nabla\vth(t_n+s,x)|^2\right)\
    dx\ ds=0.
    \end{multline*}
Letting $t_n\to+\infty$ yields
$$|\nabla\mu(t_n+\cdot,\cdot)|,|\nabla\vth(t_n+\cdot,\cdot)|\to 0\quad\text{in $L_2([0,1]\times\Om)$}.$$
This in turn yields a subsequence $(t_{n_k})$ such that
$\nabla\mu(t_{n_k}+s),\nabla\vth(t_{n_k}+s)\to 0$ in
$L_2(\Om;\R^n)$ for a.e. $s\in [0,1]$. Hence
$\nabla\vth_\infty=0$, since the gradient is a closed operator in
$L_2(\Om;\R^n)$. This in turn yields that $\vth_\infty$ is a
constant.

\noindent Furthermore the Poincar\'{e}-Wirtinger inequality
implies that
\begin{multline*}|\mu(t_{n_k}+s^*)-\mu(t_{n_l}+s^*)|_2\\\le
C_p\Big(|\nabla\mu(t_{n_k}+s^*)-\nabla\mu(t_{n_l}+s^*)|_2+\int_\Om|\Phi'(\psi(t_{n_k}+s^*))-\Phi'(\psi(t_{n_l}+s^*))|\
dx\\
+\int_\Om|\la'(\psi(t_{n_k}+s^*))\vth(t_{n_k}+s^*)-\la'(\psi(t_{n_l}+s^*))\vth(t_{n_l}+s^*)|\
dx,\end{multline*} for some $s^*\in [0,1]$. Taking the limit
$k,l\to\infty$ we see that $\mu(t_{n_k}+s^*)$ is a Cauchy sequence
in $L_2(\Om)$, hence it admits a limit, which we denote by
$\mu_\infty$. In the same manner as for $\vth_\infty$ we therefore
obtain $\nabla\mu_\infty=0$, hence $\mu_\infty$ is a constant.
Observe that the relation
$$\mu_\infty=\frac{1}{|\Om|}\left(\int_\Om(\Phi'(\psi_\infty)-\la'(\psi_\infty)\vth_\infty)\
dx\right)$$ is valid. Multiplying $\eqref{PFasym}_1$ by a function
$\ph\in H_2^1(\Om)$ and integrating by parts we obtain
\begin{multline*}
(\mu(t_{n_k}+s^*),\ph)_2=(\nabla\psi(t_{n_k}+s^*),\nabla\ph)_2+\\
(\Phi'(\psi(t_{n_k}+s^*)),\ph)_2-(\la'(\psi(t_{n_k}+s^*))\vth(t_{n_k}+s^*),\ph)_2.
\end{multline*}
 As $t_{n_k}\to\infty$ it follows that \beq\label{RegPsiInftyPF}
(\mu_\infty,\ph)_2=(\nabla\psi_\infty,\nabla\ph)_2+(\Phi'(\psi_\infty),\ph)_2-\vth_\infty(\la'(\psi_\infty),\ph)_2.
\eeq By the Lax-Milgram theorem the bounded, symmetric and
elliptic form
$$a(u,v):=\int_\Om\nabla u\nabla v\
dx,$$ defined on the space $V_1\times V_1$ induces a bounded
operator $A:V_1\to V_1^*$ with nonempty resolvent, such that
$$a(u,v)=\langle Au,v\rangle_{V_1^*,V_1},$$
for all $(u,v)\in V_1\times V_1$. It is well-known that the domain
of the part $A_p$ of the operator $A$ in
$$\X_p=\{u\in L_p(\Om):\int_\Om u\
dx=0\}$$ is given by
$$D(A_p)=\{u\in \X_p\cap H_p^2(\Om),\ \pa_\nu u=0\}.$$ Going back to
\eqref{RegPsiInftyPF} we obtain from (H1) and (H2) that
$\psi_\infty\in D(A_q)$, where $q=6/(\beta+2)$. Since $q>6/5$ we
may apply a bootstrap argument to conclude $\psi_\infty\in
D(A_2)$. Integrating \eqref{RegPsiInftyPF} by parts, assertion
2 follows.

In order to prove 3.\ , we make use of \eqref{1stder} to obtain
\begin{align*}\langle
L'(\psi_\infty,\vth_\infty)&,(h,k)\rangle_{V^*,V}\\
&=\langle
E'(\psi_\infty,\vth_\infty),(h,k)\rangle_{V^*,V}-\vth_\infty\langle
F'(\psi_\infty,\vth_\infty),(h,k)\rangle_{V^*,V}\\
&=\int_\Om(-\De\psi_\infty+\Phi'(\psi_\infty))h\
dx+\int_\Om\beta'(\vth_\infty)k\
dx\\
&\hspace{3cm}-\vth_\infty\int_\Om(\la'(\psi_\infty)h+b'(\vth_\infty)k)\
dx\\
&=\int_\Om\mu_\infty h\ dx=0,
\end{align*}
for all $(h,k)\in V$, since $\mu_\infty$ and $\vth_\infty$ are
constant. A continuity argument finally yields the last statement
of the proposition.

\epr \noindent The following result is crucial for the proof of
convergence.
\begin{prop}[Lojasiewicz-Simon inequality]\label{LS}
Let $(\psi_\infty,\vth_\infty)\in\om(\psi,\vth)$ and assume
(H1)-(H5). Then there exist constants $s\in (0,\frac{1}{2}],
C,\de>0$ such that
$$|L(u,v)-L(\psi_\infty,\vth_\infty)|^{1-s}\le C|L'(u,v)|_{V^*},$$
whenever $|(u,v)-(\psi_\infty,\vth_\infty)|_V\le \de$.
\end{prop}
\bpr We show first that $\dim
N(L''(\psi_\infty,\vth_\infty))<\infty$. By \eqref{2ndder} we obtain
\begin{align*}
\langle
L''(\psi_\infty,\vth_\infty)(h_1,k_1),&(h_2,k_2)\rangle_{V^*,V}\\
&=\int_\Om\nabla h_1\nabla h_2\
dx+\int_\Om\Phi''(\psi_\infty)h_1h_2\
dx+\int_\Om\beta''(\vth_\infty)k_1k_2\ dx\\
&-\overline{k_1}\int_\Om(\la'(\psi_\infty)h_2+b'(\vth_\infty)k_2)\
dx\\
&-\overline{k_2}\int_\Om(\la'(\psi_\infty)h_1+b'(\vth_\infty)k_1)\
dx\\
&-\overline{\vth_\infty}\int_\Om(\la''(\psi_\infty)h_1h_2+b''(\vth_\infty)k_1k_2)\
dx.
\end{align*}
Since $\beta''(s)=b'(s)+s b''(s)$ and $\vth_\infty\equiv const$ we
have
\begin{align*}
\langle
L''(\psi_\infty&,\vth_\infty)(h_1,k_1),(h_2,k_2)\rangle_{V^*,V}\\
&=\int_\Om\nabla h_1\nabla h_2\
dx+\int_\Om\left(\Phi''(\psi_\infty)h_1-\overline{k_1}\la'(\psi_\infty)-\vth_\infty\la''(\psi_\infty)h_1\right)h_2\
dx\\
&+\int_\Om(b'(\vth_\infty)(k_1-2\overline{k_1})-\overline{\la'(\psi_\infty)h_1})k_2\
dx
\end{align*}
for all $(h_j,k_j)\in V$. If $(h_1,k_1)\in
N(L''(\psi_\infty,\vth_\infty))$, it follows that
$$b'(\vth_\infty)(k_1-2\overline{k_1})-\overline{\la'(\psi_\infty)h_1}=0.$$
It is obvious that a solution $k_1$ to this equation must be
constant, hence it is given by
\beq\label{k1}k_1=-(b'(\vth_\infty))^{-1}\overline{\la'(\psi_\infty)h_1},\eeq
where we also made use of (H4). Concerning $h_1$ we have
\beq\label{eqh1}\langle
Ah_1,h_2\rangle_{V_1^*,V_1}=\int_\Om(k_1\la'(\psi_\infty)+\vth_\infty\la''(\psi_\infty)h_1-\Phi''(\psi_\infty)h_1)h_2\
dx,\eeq since $k_1$ is constant. By Proposition \ref{propomlim} it
holds that $\psi_\infty\in D(A_2)\hookrightarrow L_\infty(\Om)$,
hence $Ah_1\in H_1$, which means that $h_1\in D(A_2)$ and from
\eqref{eqh1} we obtain
$$A_2
h_1+P(\Phi''(\psi_\infty)h_1-\vth_\infty\la''(\psi_\infty)h_1-k_1\la'(\psi_\infty))=0,$$
where $P$ denotes the projection $P:H_2\to H_1$, defined by
$Pu=u-\overline{u}$. It is an easy consequence of the embedding
$D(A_2)\hookrightarrow L_\infty(\Om)$ that the linear operator
$B:H_1\to H_1$ given by
$$Bh_1=P(\Phi''(\psi_\infty)h_1-\vth_\infty\la''(\psi_\infty)h_1-k_1\la'(\psi_\infty))$$
is bounded, where $k_1$ is given by \eqref{k1}. Furthermore the
operator $A_2$ defined in the proof of Proposition \ref{propomlim}
is invertible, hence $A_2^{-1}B:H_1\to D(A_2)$ is a compact operator
by compact embedding and this in turn yields that $(I+A_2^{-1}B)$ is
a Fredholm operator. In particular it holds that $\dim
N(I+A_2^{-1}B)<\infty$, whence $N(L''(\psi_\infty,\vth_\infty))$ is
finite dimensional and moreover
$$N(L''(\psi_\infty,\vth_\infty))\subset
D(A_2)\times (H_2^r(\Om)\cap L_\infty(\Om))\hookrightarrow
L_\infty(\Om)\times L_\infty(\Om).$$ By Hypothesis
(H5), the restriction of $L'$ to the space $D(A_2)\times
(H_2^r(\Om)\cap L_\infty(\Om))$ is analytic in a neighbourhood of $(\psi_\infty,\theta_\infty)$. For the definition of analyticity in
Banach spaces we refer to \cite[Section 3]{Chi}. Now the claim
follows from \cite[Theorem 3.10 \& Corollary 3.11]{Chi}.

\epr


\noindent Let us now state the main result of this section.
\begin{thm}
Assume (H1)-(H5) and let $(\psi,\vth)$ be a global solution of \eqref{PFasym}. Then the limits
$$\lim_{t\to\infty}\psi(t)=:\psi_\infty,\quad\text{and}\quad
\lim_{t\to\infty}\vth(t)=:\vth_\infty=const$$ exist in
$H_2^1(\Om)$ and $H_2^r(\Om),\ r\in (0,1)$, respectively, and
$(\psi_\infty,\vth_\infty)$ is a strong solution of the stationary
problem \eqref{statsys}.
\end{thm}
\bpr Since by Proposition \ref{propomlim} the $\om$-limit set is
compact, we may cover it by a union of \emph{finitely} many balls
with center $(\ph_i,\th_i)\in\om(\psi,\vth)$ and radius $\de_i>0$,
$i=1,\ldots,N$. Since $L(u,v)\equiv L_\infty$ on $\om(\psi,\vth)$
and each $(\ph_i,\th_i)$ is a critical point of $L$, there are
\emph{uniform} constants $s\in(0,\frac{1}{2}]$, $C>0$ and an open
set $U\supset \om(\psi,\vth)$, such that
\beq\label{CHGLS}|L(u,v)-L_\infty|^{1-s}\le C|L'(u,v)|_{V^*},\eeq
for all $(u,v)\in U$. Define
$H:\R_+\to\R_+$ by
$$H(t):=(L(\psi(t),\vth(t))-L_\infty)^s.$$
The function $H$ is nonincreasing and $\lim_{t\to\infty}H(t)=0$, since $L(\psi(t),\vartheta(t))=E(\psi(t),\vartheta(t))$ and since $E$ is a strict Lyapunov functional for \eqref{PFasym}, which follows from \eqref{eneq}. Furthermore we have
$\lim_{t\to\infty}\dist((\psi(t),\vth(t)),\om(\psi,\vth))=0$, i.e.
there exists $t^*\ge 0$, such that $(\psi(t),\vth(t))\in U$,
for all $t\ge t^*$. Next, we compute and estimate the time
derivative of $H$. By \eqref{eneq} and Proposition \ref{LS} we
obtain
\begin{align}\label{PFA6}
-\difft\ H(t)&=s\left(-\difft\
L(\psi(t),\vth(t))\right)|L(\psi(t),\vth(t))-L_\infty|^{s-1}\nn\\
&\ge
C\frac{|\nabla\mu(t)|_2^2+|\nabla\vth(t)|_2^2}{|L'(\psi(t),\vth(t))|_{V^*}}
\end{align}
So have to estimate the term $|L'(\psi(t),\vth(t))|_{V^*}$.
For convenience we will write $\psi=\psi(t)$ and $\vth=\vth(t)$.
From \eqref{1stder} we obtain with $\bar{h}=0$
\beq\label{1}\begin{split} \langle
L'(\psi,\vth)&,(h,k)\rangle_{V^*,V}\\
&=\int_\Om(-\De\psi+\Phi'(\psi))h\
dx+\int_\Om\vth b'(\vth)k\
dx-\overline{\vth}\int_\Om(\la'(\psi)h+b'(\vth)k)\
dx\\
&=\int_\Om(\mu-\overline{\mu}) h\
dx+\int_\Om(\vth-\overline{\vth})\la'(\psi)h\
dx+\int_\Om(\vth-\overline{\vth})b'(\vth)k\ dx\end{split}\eeq An
application of the H\"{o}lder and Poincar\'{e}-Wirtinger inequality yields
the estimates
\begin{align}\label{2}
|\int_\Om(\vth-\overline{\vth})\la'(\psi)h\ dx|&\le
|\la'(\psi)|_{\infty}|\vth-\overline{\vth}|_2|h|_2\le
c|\nabla\vth|_2|h|_2,
\end{align}
\begin{align}\label{3}
|\int_\Om(\vth-\overline{\vth})b'(\vth)k\ dx|&\le
|b'(\vth)|_{\infty}|\vth-\overline{\vth}|_2|k|_2\le
c|\nabla\vth|_2|k|_2
\end{align}
and
\begin{align}\label{4}
|\int_\Om(\mu-\overline{\mu})h\ dx|&\le c|\nabla\mu|_2|h|_2,
\end{align}
whence we obtain
$$|L'(\psi(t),\vth(t))|_{V^*}\le C(|\nabla\mu(t)|_2+|\nabla\vth(t)|_2),$$
by taking the supremum over all functions $(h,k)\in V$ with norm
less than 1 in \eqref{1}-\eqref{4}. This in connection with
\eqref{PFA6} yields
$$-\difft H(t)\ge C(|\nabla\mu(t)|_2+|\nabla\vth(t)|_2),$$
hence $|\nabla\mu|,|\nabla\vth|\in L_1([t^*,\infty),L_2(\Om))$. Using
the equation $\pa_t\psi=\De\mu$ we see that $\pa_t\psi\in
L_1([t^*,\infty),H_2^1(\Om)^*)$, hence the limit
$$\lim_{t\to\infty}\psi(t)=:\psi_\infty$$ exists in $H_2^1(\Om)^*$ and even in $H_2^1(\Omega)$ thanks to Proposition \ref{relcomp}.
From equation $\eqref{PFasym}_2$ it follows that $\pa_t e\in
L_1([t^*,\infty);H_2^1(\Om)^*)$, where $e:=b(\vth)+\la(\psi)$,
i.e. the limit $\lim_{t\to\infty}e(t)$ exists in $H_2^1(\Om)^*$.
This in turn yields that the limit
$$\lim_{t\to\infty}b(\vth(t))=:b_\infty$$
exists in $L_2(\Omega)$, by relative compactness, cf.
Proposition \ref{relcomp}. By the monotonicity
assumption (H3) we obtain $\vth(t)=b^{-1}(b(\vth(t)))$ and thus
the limit of $\vth(t)$ as $t$ tends to infinity exists in
$L_2(\Om)$. From the relative compactness of the
orbit $\vth(\R_+)$ it follows that the limit
$$\lim_{t\to\infty}\vth(t)=:\vth_\infty$$
also exists in $H_2^r(\Om),\ r\in [0,1)$. Finally Proposition
\ref{propomlim} yields the last statement of the theorem.

\epr

\section{Appendix}

\emph{Proof of Proposition \ref{fixpointest}}

Let $(u,v),(\bar{u},\bar{v})\in \B_R(u^*,v^*)$. By Sobolev embedding it
holds that $u,\bar{u}$ and $v,\bar{v}$ are uniformly bounded in
$C^1(\overline{\Om})$ and $C(\overline{\Om})$, respectively.
Furthermore, we will use the following inequality, which has been
proven in \cite[Lemma 6.2.3]{Zach}. \beq\label{zachineqPF}
|f(w)-f(\bar{w})|_{H_p^s(L_p)}\le
\mu(T)(|w-\bar{w}|_{H_p^{s_0}(L_p)}+|w-\bar{w}|_{\infty,\infty}),\quad
0<s<s_0<1,\eeq valid for every $f\in C^{2-}(\R)$ and all
$w,\bar{w}\in\B_R^1(u^*)\cup\B_R^2(v^*)$. Here $\mu=\mu(T)$ denotes
a function, with the property $\mu(T)\to 0$ as $T\to 0$. The proof consists of several steps

(i) By H\"{o}lders inequality it holds that
\begin{align*}|\De&\Phi'(u)-\De\Phi'(\bar{u})|_{X(T)}\\
&\le|\De u\Phi''(u)-\De \bar{u}\Phi''(\bar{u})|_{X(T)}+
||\nabla u|^2\Phi'''(u)-|\nabla \bar{u}|^2\Phi'''(\bar{u})|_{X(T)}\\
&\le |\De u|_{rp,rp}|\Phi''(u)-\Phi''(\bar{u})|_{r'p,r'p}+|\De u-\De \bar{u}|_{rp,rp}|\Phi''(\bar{u})|_{r'p,r'p}\\
&\hspace{0.3cm}+T^{1/p}\left(|\nabla
u|_{\infty,\infty}^2|\Phi'''(u)-\Phi'''(\bar{u})|_{\infty,\infty}+|\nabla
u-\nabla \bar{u}|_{\infty,\infty}|\Phi'''(\bar{u})|_{\infty,\infty}\right)\\
&\le T^{1/r'p}\left(|\De
u|_{rp,rp}|\Phi''(u)-\Phi''(\bar{u})|_{\infty,\infty}+|\De u-\De
\bar{u}|_{rp,rp}|\Phi''(\bar{u})|_{\infty,\infty}\right)\\
&\hspace{0.3cm}+T^{1/p}\left(|\nabla
u|_{\infty,\infty}^2|\Phi'''(u)-\Phi'''(\bar{u})|_{\infty,\infty}+|\nabla
u-\nabla
\bar{u}|_{\infty,\infty}|\Phi'''(\bar{u})|_{\infty,\infty}\right),
\end{align*}
since $u,\bar{u}\in C(J;C^1(\overline{\Om}))$. We have
$$\De w\in H_p^{\th_2/2}(J;H_p^{2(1-\th_2)}(\Om))\hookrightarrow L_{r p}(J\times\Om),\quad
\th_2\in[0,1],$$ for every function $w\in E_1(T)$, since $r>1$ may
be chosen close to 1. Therefore we obtain
$$|\De\Phi'(u)-\De\Phi'(\bar{u})|_{X(T)}\le \mu(T)\left(R+|u^*|_1\right)|u-\bar{u}|_1,$$
due to the assumption $\Phi\in C^{4-}(\R)$.

(ii) Consider the term $(\la'(\psi_0)-\la'(u))\De
v-(\la'(\psi_0)-\la'(\bar{u}))\De\bar{v}$.
\begin{align*}
|(\la'(\psi_0)-\la'(u))&\De v-(\la'(\psi_0)-\la'(\bar{u}))\De\bar{v}|_{X(T)}\\
&\le |(\la'(\psi_0)-\la'(u))\De(v-\bar{v})|_{X(T)}+|(\la'(u)-\la'(\bar{u}))\De\bar{v}|_{X(T)}\\
&\le |\psi_0-u|_{\infty,\infty}|v-\bar{v}|_{E_2(T)}+|u-\bar{u}|_{\infty,\infty}|\bar{v}|_{E_2(T)}\\
&\le
(|\psi_0-u^*|_{\infty,\infty}+|u^*-u|_{\infty,\infty})|v-\bar{v}|_{E_2(T)}\\
&\hspace{1cm}+|u-\bar{u}|_{E_1(T)}(|\bar{v}-v^*|_{E_2(T)}+|v^*|_{E_2(T)})\\
&\le C(\mu(T)+R)|(u,v)-(\bar{u},\bar{v})|_1,
\end{align*}
since $\la\in C^{4-}(\R)$. Next, we consider the term
$\nabla(\la'(\psi_0)-\la'(u))\nabla
v-\nabla(\la'(\psi_0)-\la'(\bar{u}))\nabla\bar{v}$. We obtain
\begin{multline*}
|\nabla(\la'(\psi_0)-\la'(u))\nabla
v-\nabla(\la'(\psi_0)-\la'(\bar{u}))\nabla\bar{v}|_{X(T)}\\
\le
|\nabla(\la'(\psi_0)-\la'(u))|_\infty|\nabla(v-\bar{v})|_{X(T)}+
|\nabla(\la'(u)-\la'(\bar{u}))|_\infty|\nabla\bar{v}|_{X(T)}.
\end{multline*}
Since
$$\nabla(\la'(\psi_0)-\la'(u))=\nabla\psi_0(\la''(\psi_0)-\la''(u))+\la''(u)(\nabla\psi_0-\nabla
u),$$ and the same for $\nabla(\la'(u)-\la'(\bar{u}))$, we may
argue as above, to conclude
\begin{multline*}
|\nabla(\la'(\psi_0)-\la'(u))|_{\infty,\infty}|\nabla(v-\bar{v})|_{X(T)}+
|\nabla(\la'(u)-\la'(\bar{u}))|_{\infty,\infty}|\nabla\bar{v}|_{X(T)}\\
\le (\mu(T)+R)|(u,v)-(\bar{u},\bar{v})|_{1}.
\end{multline*}
Finally, we estimate the remaining part with H\"{o}lder's
inequality to the result
\begin{multline}\label{locpropPF1}
|v\De(\la'(\psi_0)-\la'(u))-\bar{v}\De(\la'(\psi_0)-\la'(\bar{u}))|_{X(T)}\\
\le|v-\bar{v}|_{\infty,\infty}|\De(\la'(\psi_0)-\la'(u))|_{X(T)}+|\bar{v}|_{r'p,r'p}|\De(\la'(u)-\la'(\bar{u}))|_{rp,rp},
\end{multline} where $1/r+1/r'=1$. For the first part, we obtain
\begin{align*}
|\De&(\la'(\psi_0)-\la'(u))|_{X(T)}\\
&\le|\De\psi_0|_p|\la''(\psi_0)-\la''(u)|_{\infty,\infty}+|\De\psi_0-\De
u
|_p|\la''(u)|_{\infty,\infty}\\
&\hspace{1cm}+|\nabla\psi_0|_{\infty,\infty}^2|\la'''(\psi_0)-\la'''(u)|_{\infty,\infty}+|\la'''(u)|_{\infty,\infty}|\nabla\psi_0-\nabla
u|_{\infty,\infty}\\
&\le C(|\psi_0-u|_{\infty,\infty}+|\nabla\psi_0-\nabla
u|_{\infty,\infty}+|\De\psi_0-\De u|_{p,p})\\
&\le C(\mu(T)+R),
\end{align*}
since $\psi_0\in H_p^2(\Om)\cap C^1(\overline{\Om})$ and $\la\in
C^{4-}(\R)$. For the second term in \eqref{locpropPF1} we obtain
\begin{align*}
|\De(\la'(u)&-\la'(\bar{u}))|_{rp,rp}\\
&\le |\De u|_{rp,rp}|\la''(u)-\la''(\bar{u})|_{\infty,\infty}+|\la''(\bar{u})|_{\infty,\infty}|\De
u-\De\bar{u}|_{rp,rp}\\
&\hspace{1cm}+|\nabla
u|_{\infty,\infty}^2|\la'''(u)-\la'''(\bar{u})|_{\infty,\infty}+|\la'''(\bar{u})|_{\infty,\infty}|\nabla
u-\nabla \bar{u}|_{\infty,\infty}\\
&\le C|u-\bar{u}|_{E_1(T)},
\end{align*}
since $u,\bar{u}\in C(J;C^1(\overline{\Om}))$ and $r>1$ can be
chosen close enough to 1, due to the fact that $\bar{v}\in
C(J;C(\overline{\Om}))$. Finally, we observe
$$|\bar{v}|_{r'p,r'p}\le |\bar{v}-v^*|_{r'p,r'p}+|v^*|_{r'p,r'p}\le
\mu(T)+R.$$

\noindent
(iii) For simplicity we set $f(u,v)=a_0\la'(\psi_0)-a(v)\la'(u)$.
Then we compute
\begin{align}\label{pro(v)}
|f(u,v)\pa_t u&-f(\bar{u},\bar{v})\pa_t \bar{u}|_{X(T)}\nn\\
&\le |\pa_t u(f(u,v)-f(\bar{u},\bar{v}))|_{X(T)}+|f(\bar{u},\bar{v})(\pa_tu-\pa_t\bar{u})|_{X(T)}\\
&\le (|\pa_t u-\pa_t u^*|_{X(T)}+|\pa_t
u^*|_{X(T)})|f(u,v)-f(\bar{u},\bar{v})|_{\infty,\infty}\nn\\
&\hspace{2cm}+|f(\bar{u},\bar{v})|_{\infty,\infty}|\pa_tu-\pa_t\bar{u}|_{X(T)}\nn\\
&\le
C(\mu_3(T)+R)|f(u,v)-f(\bar{u},\bar{v})|_{\infty,\infty}\nn\\
&\hspace{2cm}+|f(\bar{u},\bar{v})|_{\infty,\infty}|\pa_tu-\pa_t\bar{u}|_{X(T)}\nn.
\end{align}
Next we estimate
\begin{align*}
|f(u,v)&-f(\bar{u},\bar{v})|_{\infty,\infty}\\
&\le |a(v)(\la'(u)-\la'(\bar{u}))|_{\infty,\infty}+|\la'(\bar{u})(a(v)-a(\bar{v}))|_{\infty,\infty}\\
&\le |a(v)|_{\infty,\infty}|\la'(u)-\la'(\bar{u})|_{\infty,\infty}+|\la'(\bar{u})|_{\infty,\infty}|a(v)-a(\bar{v})|_{\infty,\infty}\\
&\le C(|u-\bar{u}|_{\infty,\infty}+|v-\bar{v}|_{\infty,\infty})\le
C|(u,v)-(\bar{u},\bar{v})|_1.
\end{align*}
Furthermore, we have
\begin{align*}
|f(\bar{u},\bar{v})|_{\infty,\infty}&\le
|a_0|_{\infty,\infty}|\la'(\psi_0)-\la'(\bar{u})|_{\infty,\infty}+|\la'(\bar{u})|_{\infty,\infty}|a_0-a(\bar{v})|_{\infty,\infty}\\
&\le C(|\psi_0-\bar{u}|_{\infty,\infty}+|\vth_0-\bar{v}|_{\infty,\infty})\\
&\le
C(|\psi_0-u^*|_{\infty,\infty}+|u^*-\bar{u}|_{\infty,\infty}+|\vth_0-v^*|_{\infty,\infty}+|v^*-\bar{v}|_{\infty,\infty})\\
&\le C(\mu(T)+R).
\end{align*}
The estimate of $(a_0-a(v))\Delta v-(a_0-a(\bar{v}))\Delta\bar{v}$ in $L_p(J;L_p(\Omega))$ can be carried out in a similar way.

(iv) We compute \begin{multline*}|(a(v)-a(\bar{v})f_2|_{X(T)}\le
|a(v)-a(\bar{v})|_{\infty,\infty}|f_2|_{X(T)}\le
|v-\bar{v}|_{\infty,\infty}|f_2|_{X(T)}\\\le\mu(T)|v-\bar{v}|_{E_2(T)}\le\mu(T)|(u,v)-(\bar{u},\bar{v})|_{1},
\end{multline*}
since $f_2\in X(T)$ is a fixed function, hence $|f_2|_{X(T)}\to 0$
as $T\to 0$.

(v) By trace theory, we obtain
    \begin{multline*}
    |\pa_\nu(\Phi'(u)-\Phi'(\bar{u}))|_{Y_1(T)}\\
    \le C|\Phi'(u)-\Phi'(\bar{u})|_{H_p^{1/2}(J;L_p(\Om))}+|\Phi'(u)-\Phi'(\bar{u})|_{L_p(J;H_p^2(\Om))}.
    \end{multline*}
The second norm has already been estimated in (i), so it remains
to estimate $\Phi'(u)-\Phi'(\bar{u})$ in $H_p^{1/2}(J;L_p(\Om))$.
Here we will use \eqref{zachineqPF}, to obtain
\begin{align*}|\Phi'(u)-\Phi'(\bar{u})|_{H_p^{1/2}(L_p)}&\le
\mu(T)(|u-\bar{u}|_{H_p^{s_0}(L_p)}+|u-\bar{u}|_{\infty,\infty})\\
&\le \mu(T)C|u-\bar{u}|_{E_1(T)}\le
\mu(T)C|(u,v)-(\bar{u},\bar{v})|_{1},
\end{align*}
since $s_0<1$.

(vi) We may apply (ii) and trace theory, to conclude that it
suffices to estimate
    \begin{multline*}
    (\la'(\psi_0)-\la'(u))v-(\la'(\psi_0)-\la'(\bar{u}))\bar{v}\\
    =(\la'(\psi_0)-\la'(u))(v-\bar{v})-(\la'(u)-\la'(\bar{u}))\bar{v}
    \end{multline*}
in $H_p^{1/2}(J;L_p(\Om))$. This yields
\begin{align*}
|(\la'&(\psi_0)-\la'(u))(v-\bar{v})|_{H_p^{1/2}(L_p)}\\
&\le|\la'(\psi_0)-\la'(u)|_{H_p^{1/2}(L_p)}|v-\bar{v}|_{\infty,\infty}+
|\la'(\psi_0)-\la'(u)|_{\infty,\infty}|v-\bar{v}|_{H_p^{1/2}(L_p)}\\
&\le
(|\la'(\psi_0)-\la'(u^*)|_{H_p^{1/2}(L_p)}+|\la'(u^*)-\la'(u)|_{H_p^{1/2}(L_p)})|v-\bar{v}|_{E_2(T)}\\
&\hspace{1cm}+
(|\psi_0-u^*|_{\infty,\infty}+|u^*-u|_{\infty,\infty})|v-\bar{v}|_{E_2(T)}\\
&\le
\left(|\la'(\psi_0)-\la'(u^*)|_{H_p^{1/2}(L_p)}+\mu(T)R+(\mu(T)+R)\right)|v-\bar{v}|_{E_2(T)}.
\end{align*}
Clearly
$\la'(\psi_0)-\la'(u^*)\in\, _0H_p^{1/2}(J;L_p(\Om))$, since
$\psi_0$ does not depend on $t$ and since $\la\in C^{4-}(\R)$.
Therefore it holds that
$$|\la'(\psi_0)-\la'(u^*)|_{H_p^{1/2}(L_p)}\to 0$$
as $T\to 0$. The second part $(\la'(u)-\la'(\bar{u}))\bar{v}$ can
be treated as follows.
\begin{align*}
|(\la'(u)&-\la'(\bar{u}))\bar{v}|_{H_p^{1/2}(L_p)}\\
&\le|\la'(u)-\la'(\bar{u})|_{H_p^{1/2}(L_p)}|\bar{v}|_{\infty,\infty}+|\la'(u)-\la'(\bar{u})|_{\infty,\infty}|\bar{v}|_{H_p^{1/2}(L_p)}\\
&\le C(\mu(T)+R+\mu(T))|u-\bar{u}|_{E_1(T)},
\end{align*}
where we applied again \eqref{zachineqPF}. This completes the
proof of the proposition.\vspace{0.2cm}

\noindent\emph{Proof of Proposition \ref{estpsit}}

\noindent Let $J_{\max}^\de:=[\de,T_{\max}]$ for some small $\de>0$.
Setting $A^2=\De_N^2$ with domain
$$D(A^2)=\{u\in H_p^4(\Om):\pa_\nu u=\pa_\nu\De u=0\ \text{on}\
\pa\Om\},$$ the solution $\psi(t)$ of equation $\eqref{PFGWP}_1$ may
be represented by the variation of parameters formula
\beq\label{varpar}\psi(t)=e^{-A^2t}\psi_0+\int_0^tAe^{-A^2(t-s)}\Big(\la'(\psi(s))\vth(s)-\Phi'(\psi(s))\Big)\
ds,\quad t\in J_{\max},\eeq where $e^{-A^2t}$ denotes the analytic
semigroup, generated by $-A^2=-\De_N^2$ in $L_p(\Om)$. By (H1), (H2)
and \eqref{enest2} it holds that
$$\Phi'(\psi)\in L_\infty(J_{\max};L_{q_0}(\Om))\quad\text{and}\quad
\la'(\psi)\in L_\infty(J_{\max};L_6(\Om)),$$ with $q_0=6/(\ga+2)$.
We then apply $A^{r},\ r\in (0,1)$, to \eqref{varpar} and make use
of semigroup theory to obtain \beq\label{estpsi1}\psi\in
L_\infty(J_{\max}^\de;H_{q_0}^{2r}(\Om)), \eeq valid for all $r\in
(0,1)$, since $q_0<6$. It follows from \eqref{estpsi1} that
$\psi\in L_\infty(J_{\max}^\de;L_{p_1}(\Om))$ if
$2r-3/{q_0}\ge-3/p_1$, and
$$\Phi'(\psi)\in L_\infty(J_{\max}^\de;L_{q_1}(\Om))\quad\text{as well as}\quad
\la'(\psi)\in L_\infty(J_{\max}^\de;L_{p_1}(\Om)),$$ with
$q_1=p_1/(\ga+2)$. Hence we have this time
$$\psi\in L_\infty(J_{\max}^\de;H_{q_1}^{2r}(\Om)),\quad r\in (0,1).$$
Iteratively we obtain a sequence $(p_n)_{n\in\N_0}$ such that
$$2r-\frac{3}{q_n}\ge-\frac{3}{p_{n+1}},\quad n\in\N_0$$
with $q_n=p_n/(\ga+2)$ and $p_0=6$. Thus the sequence
$(p_n)_{n\in\N_0}$ may be recursively estimated by
$$\frac{1}{p_{n+1}}\ge\frac{\ga+2}{p_n}-\frac{2r}{3},$$
for all $n\in\N_0$ and $r\in (0,1)$. From this definition it is
not difficult to obtain the following estimate for $1/p_{n+1}$.
\begin{align}\label{bootstrap}
\frac{1}{p_{n+1}}&\ge\frac{(\ga+2)^{n+1}}{p_0}-\frac{2r}{3}\sum_{k=0}^{n}(\ga+2)^k\nn\\
&=\frac{(\ga+2)^{n+1}}{p_0}-\frac{2r}{3}\left(\frac{(\ga+2)^{n+1}-1}{\ga_1+1}\right)\nn\\
&=(\ga+2)^{n+1}\left(\frac{1}{p_0}-\frac{2r}{3\ga+3}\right)+\frac{2r}{3\ga+3},\quad
n\in\N_0.
\end{align}
By the assumption (H1) on $\ga$ we see that the term in brackets
is negative if $r\in (0,1)$ is sufficiently close to 1 and
therefore, after finitely many steps the entire right side of
\eqref{bootstrap} is negative as well, whence we may choose $p_n$
arbitrarily large or we may even set $p_n=\infty$ for $n\ge N$ and
a certain $N\in\N_0$. In other words this means that for those
$r\in (0,1)$ we have \beq\label{highordest1}\psi\in
L_\infty(J_{\max}^\de;H_p^{2r}(\Om)),\eeq for all $p\in
[1,\infty]$. It is important, that we can achieve this result in
\emph{finitely} many steps!

Next we will derive an estimate for $\pa_t\psi$. For all
forthcoming calculations we will use the abbreviation
$\psi=\psi(t)$ and $\vth=\vth(t)$. Since we only have estimates on
the interval $J_{\max}^\de$, we will use the following solution
formula.
$$\psi(t)=e^{-A^2(t-\de)}\psi_\de+\int_0^{t-\de}Ae^{-A^2 s}\Big(\la'(\psi)\vth-\Phi'(\psi)\Big)(t-s)\
ds,\quad t\in J_{\max}^\de$$ where $\psi_\de:=\psi(\de)$.
Differentiating with respect to $t$, we obtain
\begin{multline}\label{psit}
\pa_t\psi(t)=
A\int_0^{t-\de}e^{-A^2s}(\la''(\psi)\vth\pa_t\psi+\la'(\psi)\pa_t\vth-\Phi''(\psi)\pa_t\psi)(t-s)\ ds\\
+F(t,\psi_\de,\vth_\de),
\end{multline}
for all $t\ge \de$ and with
$$F(t,\psi_\de,\vth_\de):=Ae^{-A^2(t-\de)}(\la'(\psi_\de)\vth_\de-\Phi'(\psi_\de))-A^2e^{-A^2 (t-\de)}\psi_\de.$$
Let us discuss the function $F$ in detail. By the trace theorem we have
$\psi_\de\in B_{pp}^{4-4/p}(\Om)$ and $\vth_{\de}\in
B_{pp}^{2-2/p}(\Om)$. Since we assume $p>(n+2)/2$, it holds that
$\psi_\de,\vth_\de\in L_\infty(\Om)$. Furthermore, the semigroup
$e^{-A^2 t}$ is analytic. Therefore there exist some constants
$C>0$ and $\om\in\R$ such that
$$|F(t,\psi_\de,\vth_\de)|_{L_p(\Om)}\le
C\left(\frac{1}{(t-\de)^{1/2}}+\frac{1}{t-\de}\right)e^{\om t},$$
for all $t>\de$. This in turn implies that
$$F(\cdot,\psi_\de,\vth_\de)\in L_p(J_{\max}^{\de'}\times\Om)$$
for all $p\in (1,\infty)$, where $0<\de<\de'<T_{\max}$. We will
now use equations $\eqref{PFasym}_{1,2}$ to rewrite the integrand
in \eqref{psit} in the following way.
\begin{align}\label{est2}
(\la''(\psi)\vth&-\Phi''(\psi))\pa_t\psi+\la'(\psi)\pa_t\vth\nn\\
&=(\la''(\psi)\vth-\Phi''(\psi))\De\mu+\frac{\la'(\psi)}{b'(\vth)}\De\vth
-\frac{\la'(\psi)^2}{b'(\vth)}\De\mu\nn\\
&=\diver\left[\left(\la''(\psi)\vth-\frac{\la'(\psi)^2}{b'(\vth)}-\Phi''(\psi)\right)\nabla\mu\right]+
\diver\left[\frac{\la'(\psi)}{b'(\vth)}\nabla\vth\right]\\
&-\nabla\left(\la''(\psi)\vth-\frac{\la'(\psi)^2}{b'(\vth)}-\Phi''(\psi)\right)\cdot\nabla\mu-\nabla\frac{\la'(\psi)}{b'(\vth)}\cdot\nabla\vth.\nn
\end{align}
Thus we obtain a decomposition of the following form
\begin{multline*}
(\la''(\psi)\vth-\Phi''(\psi))\pa_t\psi+\la'(\psi)\pa_t\vth\\
=\diver(f_\mu\nabla\mu+f_\vth\nabla\vth)+
g_\mu\nabla\mu+g_\vth\nabla\vth+h_\mu\nabla\vth\nabla\mu+h_\vth|\nabla\vth|^2,
\end{multline*}
with
\begin{align*}
f_\mu:=\la''(\psi)\vth-\frac{\la'(\psi)^2}{b'(\vth)}-\Phi''(\psi),&\quad
f_\vth:=\frac{\la'(\psi)}{b'(\vth)},\\
g_\mu:=-\left(\la'''(\psi)\vth-2\frac{\la'(\psi)\la''(\psi)}{b'(\vth)}-\Phi''(\psi)\right)\nabla\psi,&\quad
g_\vth:=-\frac{\la''(\psi)}{b'(\vth)}\nabla\psi,\\
h_\mu:=\la''(\psi)-\frac{b''(\vth)\la'(\psi)^2}{b'(\vth)^2},&\quad
h_\vth:=\frac{b''(\vth)\la'(\psi)}{b'(\vth)^2}.
\end{align*}
By Assumption (H3) and the first part of the proof it holds that
$f_j,g_j,h_j\in L_\infty(J_{\max}^\de\times\Om)$ for each
$j\in\{\mu,\vth\}$ and this in turn yields that
\begin{align*}\diver(f_\mu\nabla\mu+f_\vth\nabla\vth)&\in
L_2(J_{\max}^\de;H_2^1(\Om)^*),\\
g_\mu\cdot\nabla\mu+g_\vth\cdot\nabla\vth&\in L_2(J_{\max}^\de\times\Om),\\
h_\mu\nabla\vth\cdot\nabla\mu+h_\vth|\nabla\vth|^2&\in
L_1(J_{\max}^\de\times\Om),
\end{align*}
where we also made use of \eqref{enest2}. Setting
$$
T_1=Ae^{-A^2t}\ast \diver(f_\mu\nabla\mu+f_\vth\nabla\vth),\quad
T_2=Ae^{-A^2t}\ast (g_\mu\cdot\nabla\mu+g_\vth\cdot\nabla\vth)
$$
and
$$T_3=Ae^{-A^2t}\ast
(h_\mu\nabla\vth\cdot\nabla\mu+h_\vth|\nabla\vth|^2),$$ we may
rewrite \eqref{psit} as
$$\pa_t\psi=T_1+T_2+T_3+F(t,\psi_0,\vth_0).$$
Going back to \eqref{psit} we obtain
\begin{align*}&T_1\in
H_2^{1/4}(J_{\max}^\de;L_2(\Om))\cap L_2(J_{\max}^\de;H_2^1(\Om))\hookrightarrow L_2(J_{\max}^\de\times\Om),\\
&T_2\in H_2^{1/2}(J_{\max}^\de;L_2(\Om))\cap L_2(J_{\max}^\de;H_2^2(\Om))\hookrightarrow L_2(J_{\max}^\de\times\Om),\quad\text{and}\\
&F(\cdot,\psi_\de,\vth_\de)\in L_2(J_{\max}^{\de'}\times\Om).
\end{align*}
Observe that we do not have full regularity for $T_3$ since $A$
has no maximal regularity in $L_1(\Om)$, but nevertheless we
obtain
$$T_3\in H_1^{1/2-}(J_{\max}^\de;L_1(\Om))\cap
L_1(J_{\max}^\de;H_1^{2-}(\Om)).$$ Here we used the notation
$H_p^{s-}:=H_p^{s-\ep}$ and $\ep>0$ is sufficiently small. An
application of the mixed derivative theorem then yields
$$H_1^{1/2-}(J_{\max}^\de;L_1(\Om))\cap
L_1(J_{\max}^\de;H_1^{2-}(\Om))\hookrightarrow
L_p(J_{\max}^\de;L_2(\Om)),$$ if $p\in (1,8/7)$, whence
$$\pa_t\psi\in L_2(J_{\max}^{\de'}\times\Om)+L_p(J_{\max}^{\de'};L_2(\Om))$$
for some $1<p<8/7$. Now we go back to \eqref{est2} where we
replace this time only $\pa_t\vth$ by the differential equation
$\eqref{PFasym}_2$ to obtain
\begin{align*}\label{est3}
(\la''(\psi)\vth&-\Phi''(\psi))\pa_t\psi+\la'(\psi)\pa_t\vth\\
&=\left(\la''(\psi)\vth-\Phi''(\psi)-\frac{\la'(\psi)^2}{b'(\vth)}\right)\pa_t\psi\\
&+\diver\left[\frac{\la'(\psi)}{b'(\vth)}\nabla\vth\right]-\frac{\la''(\psi)}{b'(\vth)}\nabla\psi\cdot\nabla\vth
+\frac{\la'(\psi)b''(\vth)}{b'(\vth)^2}|\nabla\vth|^2\\
&=f\pa_t\psi+\diver\left[g\nabla\vth\right]+h\cdot\nabla\vth+k|\nabla\vth|^2.
\end{align*}
Rewrite \eqref{psit} in the following way
\beq\label{S}\pa_t\psi=S_1+S_2+S_3+S_4+F(t,\psi_0,\vth_0),\eeq
where the functions $S_j$ are defined in the same manner as $T_j$.
Since $f,g,h\in L_\infty(J_{\max}^\de\times\Om)$ it follows again
from regularity theory that
\begin{multline*}
S_1\in H_2^{1/2}(J_{\max}^{\de'};L_2(\Om))\cap
L_2(J_{\max}^{\de'};H_2^2(\Om))\\
+H_p^{1/2}(J_{\max}^{\de'};L_2(\Om))\cap L_p(J_{\max}^{\de'};H_2^2(\Om)),
\end{multline*}
    $$S_2\in H_2^{1/4}(J_{\max}^{\de'};L_2(\Om))\cap L_2(J_{\max}^{\de'};H_2^1(\Om)),$$
    $$S_3\in H_2^{1/2}(J_{\max}^{\de'};L_2(\Om))\cap L_2(J_{\max}^{\de'};H_2^2(\Om)),$$
and it can be readily verified that
$$H_p^{1/2}(J_{\max}^{\de'};L_2(\Om))\cap L_p(J_{\max}^{\de'};H_2^2(\Om))\hookrightarrow
L_2(J_{\max}^{\de'}\times\Om),$$ whenever $p\in[1,2]$. Now we turn
our attention to the term $S_4=Ae^{-A^2t}\ast k|\nabla\vth|^2$.
First we observe that by the mixed derivative theorem the
embedding
$$Z_q:=H_q^{1/2-}(J_{\max}^{\de'};L_1(\Om))\cap
L_q(J_{\max}^{\de'};H_1^{2-}(\Om))\hookrightarrow
L_2(J_{\max}^{\de'}\times\Om)$$ is valid, provided that
$q\in(8/5,2]$. Hence it holds that
$$|S_4|_{2,2}\le C|S_4|_{Z_q}\le C|k|\nabla\vth|^2|_{q,1}\le
C|\nabla\vth|_{2q,2}^2,$$ with some constant $C>0$. Taking the
norm of $\pa_t\psi$ in $L_2(J_{\max}^{\de'}\times\Om)$ we obtain
from \eqref{S}
$$|\pa_t\psi|_{2,2}\le
C\left(\sum_{j=1}^3|S_j|_{2,2}+|\nabla\vth|^2_{2q,2}+|F(\cdot,\psi_\de,\vth_\de)|_{2,2}\right).$$
The Gagliardo-Nirenberg inequality in connection with
\eqref{enest2} yields the estimate
$$|\nabla\vth|_{2q,2}^2\le
c|\nabla\vth|_{2,2}^{2a}|\nabla\vth|_{\infty,2}^{2(1-a)}\le
c|\nabla\vth|_{\infty,2}^{2(1-a)},$$ provided that $a=1/q$.
Multiply $\eqref{PFGWP}_2$ by $\pa_t\vth$ and integrate by parts
to the result
    \begin{multline*}
    \int_\Om b'(\vth(t,x))|\pa_t\vth(t,x)|^2\
    dx+\frac{1}{2}\difft|\nabla\vth(t)|_{2}^2=-\int_\Om
    \la'(\psi(t,x))\pa_t\psi(t,x)\pa_t\vth(t,x)\ dx.
    \end{multline*}
Making use of (H3) and Young's inequality we obtain
\beq\label{enest4}C_1|\pa_t\vth|_{2,2}^2+\frac{1}{2}|\nabla\vth(t)|_2^2\le
C_2(|\pa_t\psi|_{2,2}^2+|\nabla\vth_0|_2^2),\eeq after integrating
w.r.t. $t$. This in turn yields the estimate
$$|\nabla\vth|_{2q,2}^2\le
c|\nabla\vth|_{\infty,2}^{2(1-a)}\le
c(1+|\pa_t\psi|_{2,2}^{2(1-a)}).$$ In order to gain something from
this inequality we require that $2(1-a)<1$, i.e. $q$ is restricted
by $1<q<2$. Finally, if we choose $q\in (8/5,2)$ and use the
uniform boundedness of the $L_2$ norms of $S_j,\ j\in\{1,2,3\}$ we
obtain
$$|\pa_t\psi|_{2,2}\le C(1+|\pa_t\psi|_{2,2}^{2(1-a)}).$$
Since by construction $2(1-a)<1$, it follows that the $L_2$-norm
of $\pa_t\psi$ is bounded on $J_{\max}^{\de'}\times\Om$. In
particular, this yields the statement for $\vth$ by equation
\eqref{enest4}.

Now we go back to \eqref{psit} with $\de$ replaced by $\de'$. By
Assumption (H5), by the bounds $\pa_t\vth,\pa_t\psi\in
L_2(J_{\max}^{\de'};L_2(\Om))$ and by the first part of the proof
we obtain
$$\la''(\psi)\vth\pa_t\psi+\la'(\psi)\pa_t\vth-\Phi''(\psi)\pa_t\psi\in
L_2(J_{\max}^{\de'};L_2(\Om)).$$ Since the operator $A^2=\De^2$
with domain
$$D(A^2)=\{u\in H_p^4(\Om):\pa_\nu u=\pa_\nu\De u=0\}$$
has the property of maximal $L_p$-regularity (cf. \cite[Theorem 2.1]{DHP07}), we obtain from \eqref{psit}
$$\pa_t\psi-F(\cdot,\psi_{\de'},\vth_{\de'})\in
H_2^{1/2}(J_{\max}^{\de'};L_2(\Om))\cap
L_2(J_{\max}^{\de'};H_2^2(\Om))\hookrightarrow
L_r(J_{\max}^{\de'};L_r(\Om)),$$ and the last embedding is valid
for all $r\le 2(n+4)/n$. By the properties of the function $F$ it
follows
$$\pa_t\psi\in L_r(J_{\max}^{\de''};L_r(\Om)),$$
for all $r\le 2(n+4)/n$ and some $0<\de''<T_{\max}$. To obtain an
estimate for the whole interval $J_{\max}$, we use the fact that
we already have a local strong solution, i.e. $\pa_t\psi\in
L_p(0,\de'';L_p(\Om))$, $p>(n+2)/2$. The proof is complete.

\bibliographystyle{plain}
\bibliography{PenFife}

%

\end{document}